\documentclass[times,final]{elsarticle}

\usepackage{framed,multirow}

\usepackage{amssymb}
\usepackage{latexsym}
\usepackage{mathdots}

\usepackage{url}
\usepackage{xcolor}
\definecolor{newcolor}{rgb}{.8,.349,.1}
 \usepackage{blkarray,stackengine,tikz,trimclip}
\usepackage{graphicx}
\usepackage{amsmath,latexsym}
\usepackage{amsfonts}
\usepackage{bbm}
\usepackage{xcolor}
\usepackage{enumitem}

\usepackage{graphicx}
\usepackage{ulem}

\definecolor{redorange}{HTML}{F26035}
\newcommand{\jr}[1]{{\color{redorange} #1}}
\newcommand{\ignore}[1]{}
\usepackage{comment}
\usepackage{todonotes}
\usepackage[margin=1in]{geometry}
\usepackage{comment}
\usepackage{todonotes}
\newlist{todolist}{itemize}{2}
\setlist[todolist]{label=$\square$}

\begin{document}


\begin{frontmatter}

\title{Denoising Particle-In-Cell Data via Smoothness-Increasing
Accuracy-Conserving Filters with Application to Bohm Speed Computation}

\author[1]{Matthew J. Picklo}


\author[2]{Qi Tang}
\author[2]{Yanzeng Zhang}
\author[1,3]{Jennifer K. Ryan\corref{cor1}}
\cortext[cor1]{Corresponding author}
\ead{jryan@kth.se}
\author[2]{Xian-Zhu Tang}
\address[1]{Colorado School of Mines, 1500 Illinois St, Golden CO 80401, USA}
\address[2]{Theoretical Division, Los Alamos National Laboratory,
Los Alamos, NM 87545, USA}
\address[3]{KTH Royal University, Stockholm, Sweden}


\begin{abstract}
The simulation of plasma physics is computationally expensive because the underlying physical system is of high dimensions, requiring three spatial dimensions and three velocity dimensions.
One popular numerical approach is Particle-In-Cell (PIC) methods owing to its ease of implementation and favorable scalability 
in high-dimensional problems. An unfortunate drawback of the method is the introduction of statistical noise resulting from the 
use of finitely many particles. In this paper we examine the application of the Smoothness-Increasing Accuracy-Conserving (SIAC) 
family of convolution kernel filters as denoisers for moment data arising from PIC simulations. We show that SIAC filtering is a promising tool to 
denoise PIC data in the physical space as well as capture the appropriate scales in the Fourier space. 
Furthermore, we demonstrate how the application of the SIAC technique reduces the amount of information necessary 
in the computation of quantities of interest in plasma physics such as the Bohm speed.
\end{abstract}


\end{frontmatter}



 
\section{Introduction}\label{section1}

\begin{figure}[tp!]
      \centering
     \includegraphics[width=0.45\linewidth]{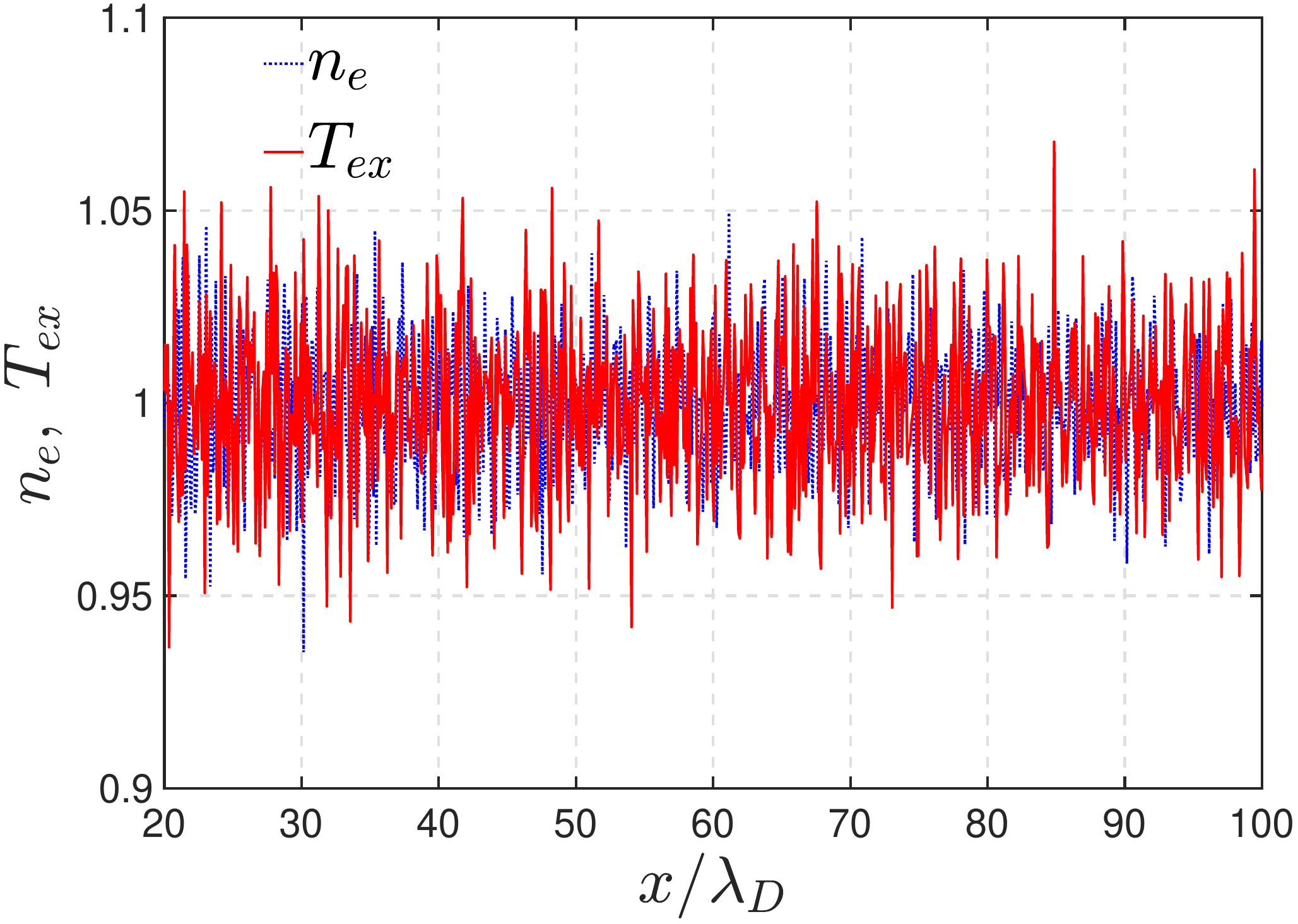} \qquad
     \includegraphics[width=0.45\linewidth]{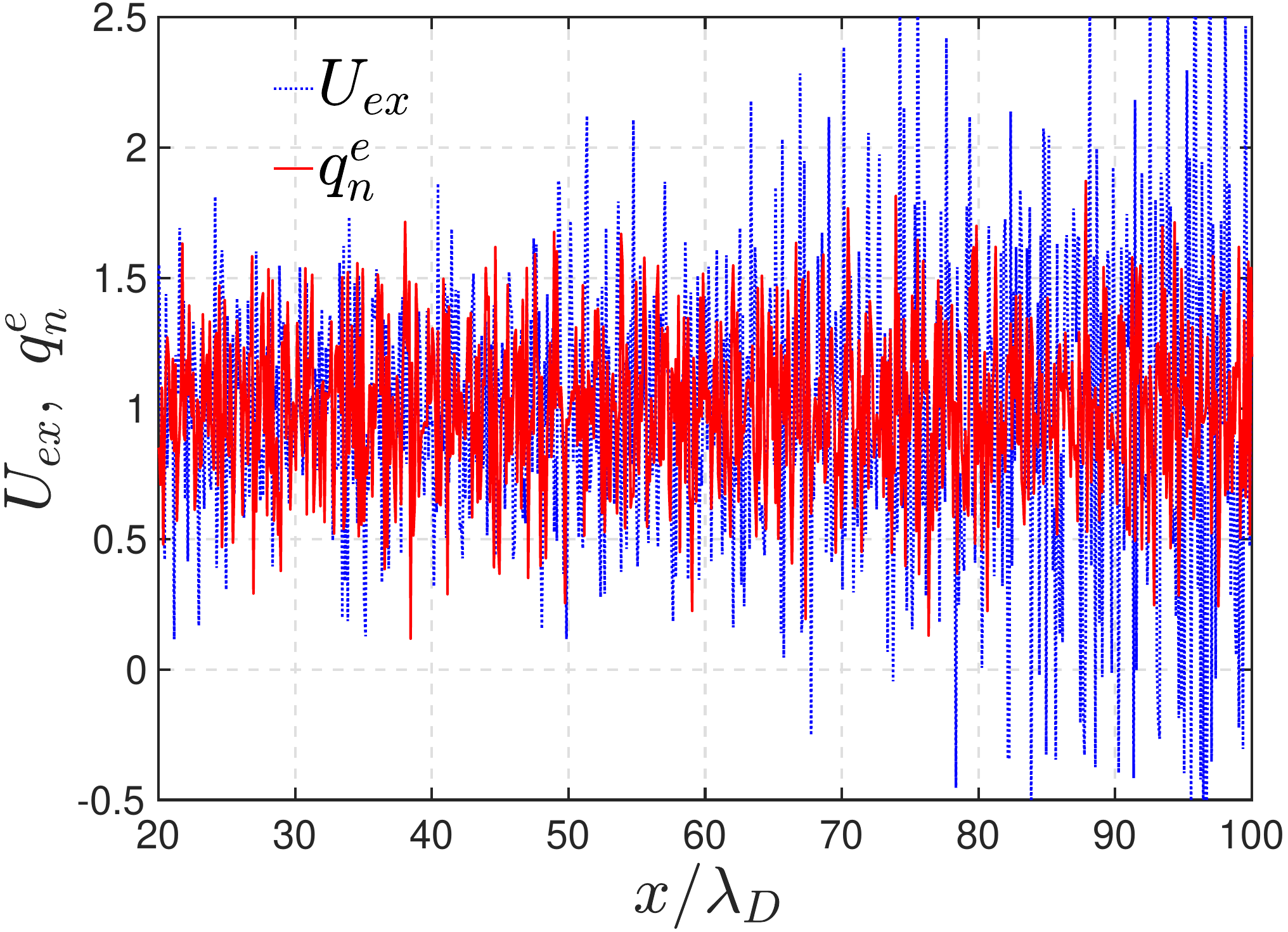}
          \caption{Electron density $n_e$ and temperature $T_{ex}$ (left figure), and the flow $U_{ex}$ and thermal conduction heat flux $q_n^e$ (right figure) from 1D3V PIC simulations of a plasma thermal quench with absorbing boundary from Ref.~\cite{zhang2023cooling}, where 5,000 particles per cell were used. All the variables are normalized by their denoised post-processed quantities. }
     \label{fig:PIC-noise}
      \end{figure}

With the high-dimensionality of the equations governing plasma
kinetics in six dimensional phase space~\cite{Krall-Trivelpiece-1973},
traditional numerical techniques for solving partial differential
equations (PDE) suffer from the curse of dimensionality.  This
motivated the development and application of Particle-In-Cell (PIC) methods
that have been a staple in plasma kinetic
simulations ~\cite{Hockney-Eastwood-1988,Birdsall-Langdon-1985}. The
particle-based approach provides simplicity and geometric flexibility,
but these benefits come at the expense of siginificant amount of
noise, which scales as $1/\sqrt{N}$ with $N$ the number of particle
markers in a cell. As an example, for a resolution of $10^{4}$ particle
markers per cell, the noise level is still about 1\%. Remarkably,
finite amounts of noise, even at substantial level, does not generally
prevent a physics-wise meaningful PIC simulation.  For example,
demonstration of Laudau damping can be achieved with a rather modest
$N$, except that the post-run physics analysis is hampered by the noisy
diagnostics.  The difficulty is aggravated by the fact that the
underlying physics is usually described, in the standard statistical
physics approach, by velocity moments of the particle distribution
function $f(\mathbf{x}, \mathbf{v}, t)$,
\begin{align}
\left< v_i^k \right> = \int v_i^k f(\mathbf{x}, \mathbf{v}, t) \; d^3\mathbf{v},
\end{align}
with $v_i$ the cartesian component of $\mathbf{v}.$
For example, the thermodynamic state variables have the density
$n(\mathbf{x},t) \equiv \left<\mathbf{v}^0\right>,$ the flow $\mathbf{U}(\mathbf{x},t) \equiv \left<
\mathbf{v}\right>/n,$ and the temperature $T(\mathbf{x},t) \equiv (m/2) \left(\left< |\mathbf{v}|^2
\right>/n - |\mathbf{U}|^2\right).$ Closure quantities such as the plasma heat flux involve
even higher order velocity moments, for example,
\begin{align}
\mathbf{q} \equiv \frac{m}{2}\left< \left(\mathbf{v} - \mathbf{U}\right)^3 \right>.
\end{align}
The general rule of thumb is that the particle noise becomes further
exaggerated in higher order moment quantities. In
Fig.~\ref{fig:PIC-noise}, we show the fluctuations in commonly used
physical quantities that are obtained by post-processing particle data
from a PIC simulation.
This high level of noise in the post-processed physical quantities
hampers our ability to gain deep insights into the underlying physics.

A recent work on the Bohm criterion analysis~\cite{li-Bohm-PRL} that
establishes the lower bound of the plasma outflow speed, also known as
the Bohm speed~\cite{Bohm}, in a wall-bounded plasma with finite
collisionality, illustrates the extent one has to go to overcome the
PIC noise in order to decipher the underlying physics.  Specifically,
the Bohm speed is predicted by theoretical analysis to take the
form of 
\begin{equation}\label{Bohm_eqn}
u_{Bohm}\equiv \sqrt{\frac{Z\beta T^{se}_{ex}+3T^{se}_{ix}}{m_i}},
\end{equation}
where 
\begin{equation}\label{eq-beta}
\beta\equiv
\frac{3-\frac{3+2\alpha}{Ze\Gamma^{se}_i}\Big(\frac{\partial
    q^i_n}{\partial
    \phi}+\frac{Q_{ii}}{E}\Big)+\frac{\alpha}{e\Gamma^{se}_e}\Big(\frac{\partial
    q^e_n}{\partial
    \phi}+\frac{Q_{ee}+Q_{ei}}{E}\Big)}{1+\frac{1+\alpha}{e\Gamma^{se}_e}\Big(\frac{\partial
    q^e_n}{\partial \phi}+\frac{Q_{ee}+Q_{ei}}{E}\Big)},
\end{equation}
with $\Gamma_{e,i}=n_{e,i}u_{e,ix}$, $\partial q^{e,i}_n/\partial
\phi=-(\partial q^{e,i}_n/\partial x)/E$, $\alpha
=-R_T/(n_e\partial T_{ex}/\partial x)$ being the thermal force coefficient, $Z$ the ion charge state, $m_i$ the ion mass, and $e$ the elementary charge.  The plasma transport physics are in the
form of local electric field $E$, which is proportional to spatial
gradient of the potential ($\phi$), particle and heat fluxes 
(which are denoted by $\Gamma_{e,i}$ and $q_n^{e,i}$, respectively),
collisional energy exchange between electrons $(Q_{ee})$ and between
electrons and ions $(Q_{ei})$, electron and ion temperatures $(T_{ex},T_{ix})$, respectively, and the thermal force ($R_T$). Here the superscript $se$
denotes the spatial location at which all quantities {{are}} to be evaluated
to find the local Bohm speed, and the subscripts $e,i$ stands for the electron and ion quantities, respectively, and $x$ represents the x-component of the variable. Most of
these quantities are simply velocity moments of the computed distribution
function, or the velocity moments of the Coulomb collisional integral,
and quantities that are derived from them. The PIC noise thus makes it
an extremely difficult task to directly compare the numerical
simulation results with theoretical predictions.
Ref.~\cite{li-Bohm-PRL} overcame this hurdle by performing long-time
simulations after a steady-state is reached, and then applying
time averaging over an extremely long-period of PIC data to control the PIC
noise to {{be a}} sufficiently low level so {{that}} a definite comparison with
Eqs.~\eqref{Bohm_eqn} and \eqref{eq-beta} {{can be done}}. This effectively increases $N$
from $10^{3-4}$ markers per cell to tens of millions per cell by
accumulating data over time.  Powerful and computationally expensive
as this approach is, many applications do not even allow a
steady-state solution. To decipher subtle transport physics from PIC
simulations, one must therefore explore alternative approaches to
effectively deal with PIC noise when modestly large $N$ and limited
time-window averaging are available.
We note that the regression and denoising approaches applied to the velocity space of PIC simulations are available, 
such as those based on Gaussian mixtures~\cite{nguyen2020adaptive, chen2021unsupervised}. Other denoising approaches for PIC code include the $\delta f$ approach~\cite{dimits1993partially,parker1993fully}, sparse grid techniques~\cite{ricketson2016sparse,muralikrishnan2021sparse}, phase space smoothing~\cite{denavit1972numerical,wang2011particle}, wavelet denoising~\cite{gassama2007wavelet}, etc.  
Our work aims to propose and explore a distinctive approach tailored for PIC moments.

Here we address this issue through the use of Smoothness-Increasing
Accuracy-Conserving (SIAC) filters which were originally developed as
a post-processor for extracting superconvergence from finite element method (FEM) based
numerical approximations \cite{Bramble1977}. The purpose of this paper is to demonstrate the effectiveness of SIAC filtering in  ``offline'' noise reduction of pointwise data arising from PIC simulations. We
will show that the application of specially tuned SIAC convolution
kernels can selectively remove high frequency oscillations without
polluting the intended low-frequency information, and thereby improve
approximation quality. Extensions of the filtering procedure for
non-periodic data are constructed, and results are shown for how SIAC
filtering schemes can reduce the quantity of data needed in computing
the Bohm speed, an important parameter in plasma physics bounding from
below the ion velocity to a charged boundary.  In these later examples, the spatial tuning requires the kernel scaling to be adjusted for both the variable being filtered as well as the location within the domain.

The rest of the paper is organized in as follows. In Sec.~\ref{Intialization} we describe a
procedure for initializing an approximation in a continuous variable
from discrete data using Lagrange polynomials. In Sec.~\ref{SIAC
  Filters} we describe the construction of symmetric and
position-dependent SIAC kernels which we convolve with our
initialization dampen high-frequency oscillations. Boundary treatments
as well as adaptive kernel scalings are also discussed. Lastly, in
Sec.~\ref{Numerical Results} we provide numerical examples
detailing the damping effects of symmetric kernel on periodic data,
adaptively scaled and position-dependent generalized spline kernels
for non-periodic data, and data reduction enabled by SIAC filtering in
Bohm speed computation.
  \begin{figure}[tp!]
  \centering
    
\begin{tabular}{c c c}

\begin{tikzpicture}[scale=2.2]
 \usetikzlibrary{calc}

 \coordinate (Origin) at (0,0);

\draw[black,fill] (-1.5,0) circle (0.5pt) node[below]{$\textbf{x}_{j-3}$};

\draw[black,fill] (-1.0,0) circle (0.5pt) node[below]{$\textbf{x}_{j-2}$};

\draw[black,fill] (-0.5,0) circle (0.5pt) node[below]{$\textbf{x}_{j-1}$};
\draw[black,fill] (0,0) circle (0.5pt) node[below]{$\textbf{x}_j$};
\draw[black,fill] (0.5,0) circle (0.5pt) node[below]{$\textbf{x}_{j+1}$};
\draw[black,fill] (1,0) circle (0.5pt) node[below]{$\textbf{x}_{j+2}$};
\draw[black,fill] (1.5,0) circle (0.5pt) node[below]{$\textbf{x}_{j+3}$};

\draw[black] (-1.75,0)--(1.75,0);

\draw[red,thick] (-1.75,0.5+-0.625)--(-1.75,0.5+-0.375);
\draw[red,thick] (-1.25,0.5+-0.625)--(-1.25,0.5+-0.375);
\draw[red,thick] (-0.75,0.5+-0.625)--(-0.75,0.5+-0.375);
\draw[red,thick] (-0.25,0.5+-0.625)--(-0.25,0.5+-0.375);
\draw[red,thick] (0.25,0.5+-0.625)--(0.25,0.5+-0.375);
\draw[red,thick] (0.75,0.5+-0.625)--(0.75,0.5+-0.375);
\draw[red,thick] (1.25,0.5+-0.625)--(1.25,0.5+-0.375);

\draw[red,thick] (1.75,0.5+-0.625)--(1.75,0.5+-0.375);

\draw[blue] (-1.75,0.5-0.125)--(-1.25,0.5+0.125);
\draw[blue] (-1.25,0.75)--(-0.75,0.75);

\draw[blue] (-0.75,0.75+0.1*0.5)--(-0.25,0.75-0.1*0.5);
\draw[blue] (-0.25,0.65-0.15*0.5)--(0.25,0.65+0.15*0.5);
\draw[blue] (0.25,0.8-0.1*0.5)--(0.75,0.8+0.1*0.5);

\draw[blue] (0.75,0.9+0.1*0.5)--(1.25,0.9-0.1*0.5);
\draw[blue] (1.25,0.8)--(1.75,0.8);

\draw[black,fill] (-1.5,0.5) circle (0.5pt) node[below]{$\textbf{f}_{j-3}$};

\draw[black,fill] (-1.0,0.75) circle (0.5pt) node[below]{$\textbf{f}_{j-2}$};

\draw[black,fill] (-0.5,0.75) circle (0.5pt) node[below]{$\textbf{f}_{j-1}$};

\draw[black,fill] (0,0.65) circle (0.5pt) node[below]{$\textbf{f}_j$};

\draw[black,fill] (0.5,0.8) circle (0.5pt) node[below]{$\textbf{f}_{j+1}$};

\draw[black,fill] (1,0.9) circle (0.5pt) node[below]{$\textbf{f}_{j+2}$};

\draw[black,fill] (1.5,0.8) circle (0.5pt) node[below]{$\textbf{f}_{j+3}$};

\end{tikzpicture}

&

&

\begin{tikzpicture}[scale=1.75]
 \usetikzlibrary{calc}

 \coordinate (Origin) at (0,0);

\draw[black,fill] (-1.5,0) circle (0.5pt) node[below]{$\textbf{x}_{j-3}$};

\draw[black,fill] (-1.0,0) circle (0.5pt) node[below]{$\textbf{x}_{j-2}$};

\draw[black,fill] (-0.5,0) circle (0.5pt) node[below]{$\textbf{x}_{j-1}$};
\draw[black,fill] (0,0) circle (0.5pt) node[below]{$\textbf{x}_j$};
\draw[black,fill] (0.5,0) circle (0.5pt) node[below]{$\textbf{x}_{j+1}$};
\draw[black,fill] (1,0) circle (0.5pt) node[below]{$\textbf{x}_{j+2}$};
\draw[black,fill] (1.5,0) circle (0.5pt) node[below]{$\textbf{x}_{j+3}$};

\draw[black,fill] (2.0,0) circle (0.5pt) node[below]{$\textbf{x}_{j+4}$};

\draw[black] (-1.75,0)--(2.25,0);

\draw[red,thick] (-1.75,0.5+-0.625)--(-1.75,0.5+-0.375);

\draw[red,thick] (-0.75,0.5+-0.625)--(-0.75,0.5+-0.375);

\draw[red,thick] (0.25,0.5+-0.625)--(0.25,0.5+-0.375);

\draw[red,thick] (1.25,0.5+-0.625)--(1.25,0.5+-0.375);

\draw[red,thick] (2.25,0.5+-0.625)--(2.25,0.5+-0.375);

\draw[blue] (-1.75,0.5-0.125)--(-0.75,0.75+0.125);

\draw[blue] (-0.75,0.75+0.1*0.5)--(0.25,0.65-0.1*0.5);

\draw[blue] (0.25,0.8-0.1*0.5)--(1.25,0.9+0.1*0.5);

\draw[blue] (1.25,0.8+0.5*0.1)--(2.25,0.7-0.1*0.5);

\draw[black,fill] (-1.5,0.5) circle (0.5pt) node[below]{$\textbf{f}_{j-3}$};

\draw[black,fill] (-1.0,0.75) circle (0.5pt) node[below]{$\textbf{f}_{j-2}$};

\draw[black,fill] (-0.5,0.75) circle (0.5pt) node[below]{$\textbf{f}_{j-1}$};

\draw[black,fill] (0,0.65) circle (0.5pt) node[below]{$\textbf{f}_j$};

\draw[black,fill] (0.5,0.8) circle (0.5pt) node[below]{$\textbf{f}_{j+1}$};

\draw[black,fill] (1,0.9) circle (0.5pt) node[below]{$\textbf{f}_{j+2}$};

\draw[black,fill] (1.5,0.8) circle (0.5pt) node[below]{$\textbf{f}_{j+3}$};

\draw[black,fill] (2.0,0.7) circle (0.5pt) node[below]{$\textbf{f}_{j+4}$};

\end{tikzpicture}

\\
(a)&& (b)

  \end{tabular}

  \caption{Example of partitioning pointwise PIC data into elements, and constructing an interpolant of the pointwise data $(f_j=f(x_j))$ on each element. The red partitions represent the superimposed mesh to the pointwise grid. The blue lines represent the piecewise linear interpolant reconstructed from a stencil containing two grid values. In (a) the elements $\tau_j = I_{j}$ consist of a single cell, and  the stencil uses information from the right adjacent cell ($\ell=0$ and $r=1$). In (b)  the elements $\tau_j = I_{2j-1} \cup I_{2j}$ are formed of a union of two neighboring cells and the stencil for each element only consists of those cells inside that element. }

      \label{fig:PartitionEX}
\end{figure}
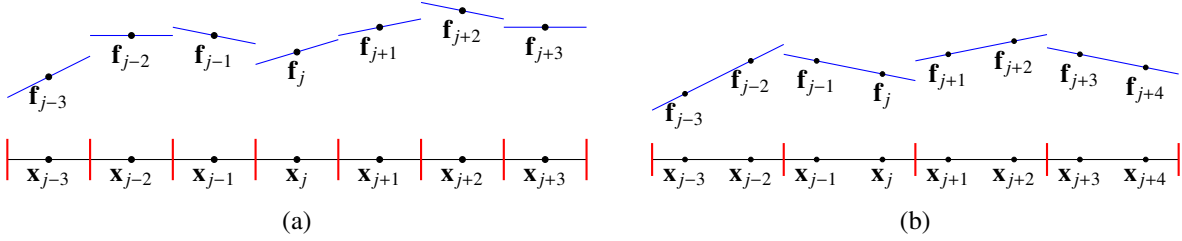

\section{Data Initialization}\label{Intialization}

Unlike traditional FEM data which is defined over the entirety of a given mesh, the data considered here arises from moments of PIC simulations and is discrete in nature. SIAC filtering makes use of a continuous convolution, therefore to apply SIAC filtering a global representation must be constructed. The methodology presented here is to superimpose a mesh over the pointwise grid where the data is defined. For simplicity, within each element of this mesh, we construct a Lagrange interpolant from the pointwise data using a stencil of a chosen size. For ease of illustration we consider the one-dimensional case.

\subsection{Constructing the mesh}

Consider a domain $\Omega=[a,b]$, and assume that the pointwise data $\{f_{j}\}^N_{j=1}$ is given on the cell centers of a grid 
\[ a<x_{1}<x_{2}<\hdots<x_{N-1}<x_{N}<b.
\]
Thinking about the data from a finite volume perspective,  introduce cells $I_j=[x_{j-1/2},x_{j+1/2}]$, $j=1,\hdots,N$, where $x_{j-1/2}=\frac{1}{2}(x_j-x_{j-1})$,  with $x_{1/2}=a$ and $x_{N+1/2}=b$. Denote the cell widths by $\Delta x_j=x_{j+1/2}-x_{j-1/2}$.
 It is possible to choose the elements of the superimposed mesh, ${\mathcal{T} =  }\{\tau_k\}_{k=1}^{N_{\tau}}$, to just be these cells, in which case  the mesh will be similar to that depicted in Fig.~\ref{fig:PartitionEX}(a). Alternatively,  the elements {can be selected} to be a union of adjacent cells, an example of which is given in Fig.~\ref{fig:PartitionEX}(b). This choice will cause no difficulty in what follows, it only alters the number of polynomial interpolants needed in the construct{ion}, $N_{\tau}$, and the minimum stencil size of these interpolants. In  the numerical implementation,  only  the cases of elements consisting of one or two cells {is considered}. 

\subsection{Constructing a piecewise polynomial interpolant}
 To convert  the discrete data into a continuous approximation on a given element $\tau{_k}$, we choose nonnegative integers $\ell$ and $r$ and introduce the stencil about element $I_j:$
\[S(j)=\{I_{j-\ell},\hdots,I_j,\hdots, I_{j+r}\}.\]
Assume here that the center cell of the stencil, $I_j$, satisfies $I_j \subset \tau{_k}$. The pointwise data contained within the cells of this stencil will be used to construct the interpolant over $\tau{_k}$. Note that the lack of assumptions about boundary conditions necessitates that  the stencil {is reduced} in size near domain boundaries {in order} to not exceed the boundaries of our domain. Setting $p=\ell+r$,  the Lagrange basis specific to  the element  is $\{L^j_q(x)\}_{q=0}^{p}$, with
\[L^j_q(x)=\prod_{\substack{n=0 \\ n\neq q}}^{p}\frac{x-x_{n+j-\ell}}{x_{q+j-\ell}-x_{n+j-\ell}}.\]
The element-specific interpolant is then given by
\[u_h(x)\Big{|}_{\tau}=\sum_{q=0}^pf_{q+j-\ell}L^j_q(x).\]
Iterating through the elements, and restricting the stencils where appropriate,   a piecewise polynomial representation of the pointwise data over the whole domain {is obtained}. In the work to follow  we find there is no significant difference in using larger stencils for the initialization procedure. To that end,  here a finite-volume perspective {is used} and  the discrete data {is treated} as cell-center values on a finite-volume style mesh. This corresponds to a $p=0$ piecewise-constant initialization.

\section{SIAC Filters: Kernel Construction and Boundary Treatments}\label{SIAC Filters}
This section  describes a convolution denoising procedure  and the construction of the SIAC kernel used in the convolution. 
For ease of illustration we again consider the one-dimensional case.

 The filtered data, $u^{\star}_h$, is obtained via convolution of the interpolant, $u_h$, with a compactly supported scaled kernel function $K_H$:
\begin{equation} \label{conv}
    u^{\star}_h(x)=K_H\star u_h\equiv\int_{\mathbb{R}}K_H(x-y) \, u_h(y)\;dy.
\end{equation} 
Here,  a family of kernel functions known as SIAC kernels is considered. These SIAC kernels are piecewise polynomial functions of compact support that satisfy consistency and moment conditions, which will be discussed  in Sec.~\ref{subsec:coeff}.  It is important to note that the compact support of the kernel function reduces the integral in Eq.~\eqref{conv} from the entire range of $\mathbb{R}$ to 
 a shifted support of the kernel. In applying the filter to non-periodic data, it is exactly this restriction  of the kernel's support that enables shifted or position-dependent kernels to constrict the convolution to within the domain of the data. 
In the following, a review of several types of SIAC kernels for applications to both periodic and non-periodic data is given. In the latter case,  a description of the position-dependent kernels developed in \cite{Ryan2003} and the generalized spline boundary kernels introduced in \cite{Ryan2015} {is given}, as well as an introduction to a novel adaptive kernel scaling 
{based on \cite{Jallepalli2019}}.

\subsection{Kernel formulation}

The knot matrix formulation of the SIAC kernel introduced in \cite{Ryan2015} is used. 
This description encompasses both the symmetric formulation for periodic data and 
the position-dependent formulation for non-periodic data. The SIAC kernel composed of $(r+1)$ B-splines of order $\ell$
 is defined via an $(r+1)\times(\ell)$ knot matrix $\mathbf{T}$ which will be given  in Sec.~\ref{subsec:knotmatrix}. 
 The composite kernel is given by
\begin{equation*} 
K_{\mathbf{T}}(x)=\sum_{\gamma=0}^{r}c_{\gamma}B^{\ell}_{\mathbf{T}_{\gamma}}(x),
\end{equation*}
where $\mathbf{T}_{\gamma}=\mathbf{T}(\gamma,:)$ denotes the $\gamma$-th row of $\mathbf{T}$ for $\gamma=0,1,\hdots,r$. 
These rows contain the knot sequences for each B-spline composing the kernel. As detailed in \cite{deBoor} the $j$-th B-spline of order $\ell$ is defined via its knot sequences $\mathbf{t}$ by the recursion relation:
\begin{equation*} 
B^{\ell}_{j}(x)=w_{j,\ell}B^{\ell-1}_{j}(x)+(1-w_{j+1,\ell})B^{\ell-1}_{j+1}(x),
\end{equation*}
where
\begin{equation*} 
w_{j,\ell}(x)=\frac{x-t_j}{t_{j+\ell-1}-t_j}.
\end{equation*}
This recursion ends with a simple characteristic function,
\begin{equation*} 
 B_j^{1}(x)=\begin{cases}
    1,&t_j\leq x< t_{j+1}\\
    0,&\text{Else.}
    \end{cases}
\end{equation*}

For the purposes of the B-splines in  the kernel,   define $B^{\ell}_{\mathbf{T}_{\gamma}}$ to be the first B-spline of order $\ell$ corresponding to knot sequence $\mathbf{T}_{\gamma}$. The typical knot matrices  considered here, with the exception of generalized spline knot matrices \cite{Ryan2015}, are of the form 
\begin{equation}\label{knot matrix}
\mathbf{T}(i,j;x^{\star})=-\frac{1}{2}(r+\ell-1)+i+j+\lambda(x^{\star}),\;\;\ i=0,\hdots,r,\;\;\; j=0,\hdots,\ell.
\end{equation}
where $\lambda (x^{\star})$ is a mesh dependent shifting function that restricts the kernel support to the computational domain. 
Note that here and in what follows,  the dependence of $\mathbf{T}_{\gamma}$ on $x^{\star}$ is suppressed in order to ease presenation.

 When employing a scaled kernel,  a multiplication of the knot matrix by a scaling function $H$ is done. 
 $H$  may depend on the location of the filtering point $x^{\star}$. Note  that, 
 occasionally,  use of the notation $K^{(r+1,\ell)}_H$ is made, which refers to the kernel function
  composed of $(r+1)$ B-splines of order $\ell$ with scaling $H$. 

\subsection{Shifting function}
\begin{figure}[tp!]
    \centering
    \includegraphics[width=0.5\linewidth]{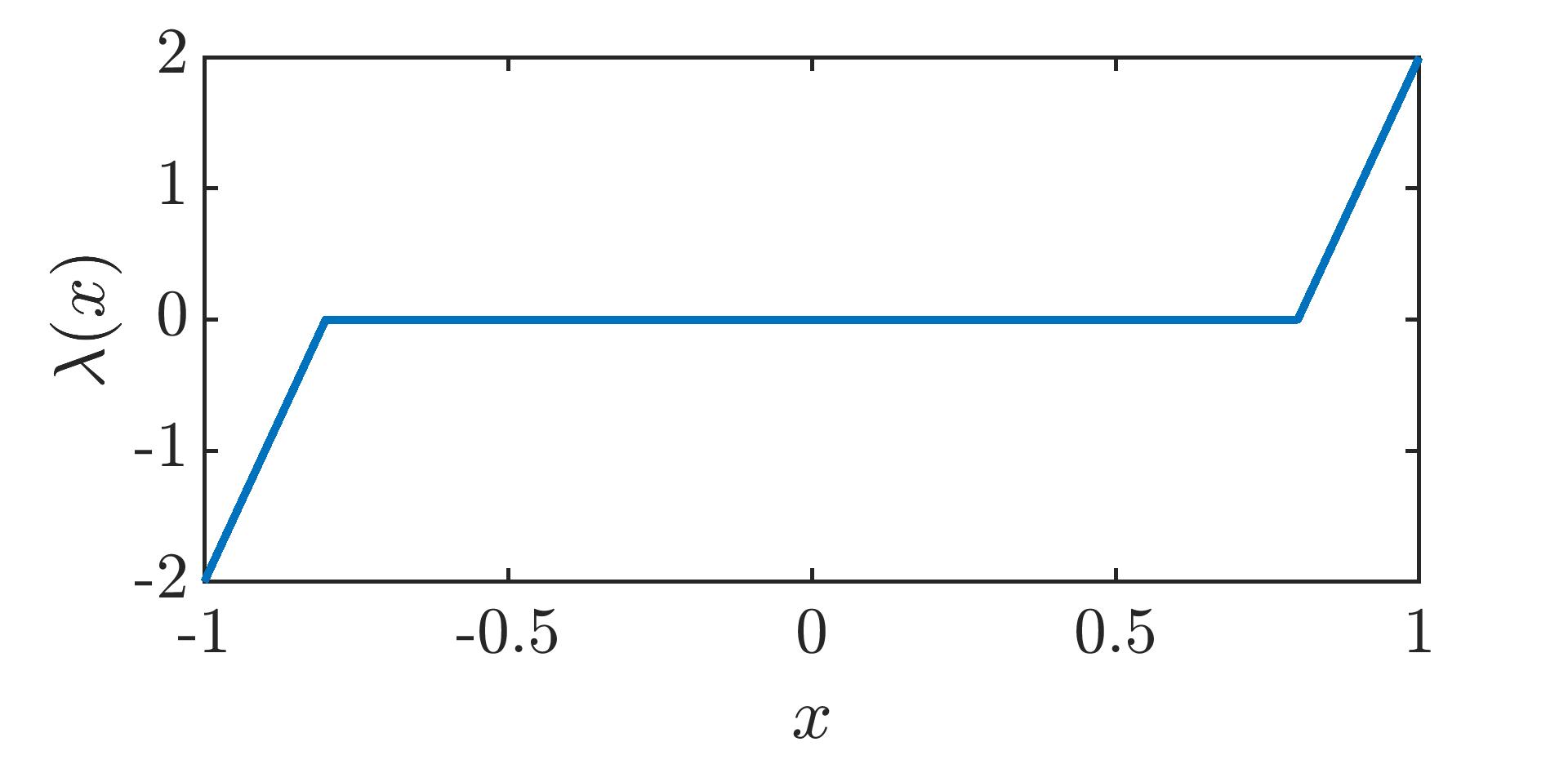}
    \caption{Demonstration of a shift function $\lambda(x^\star)$  for $K_H^{(3,2)}$ with $H(x^{\star})=0.1$ on the domain $\Omega=[-1,1]$.}
    \label{fig:shift}
\end{figure}

\begin{figure}[tp!]
    \centering
\begin{tabular}{|c|c|c|}\hline
\multicolumn{1}{|c|}{$\lambda=-2$}&\multicolumn{1}{|c|}{$\lambda=0$}&\multicolumn{1}{|c|}{$\lambda=2$}\\ \hline 
       \includegraphics[width=0.33\linewidth]{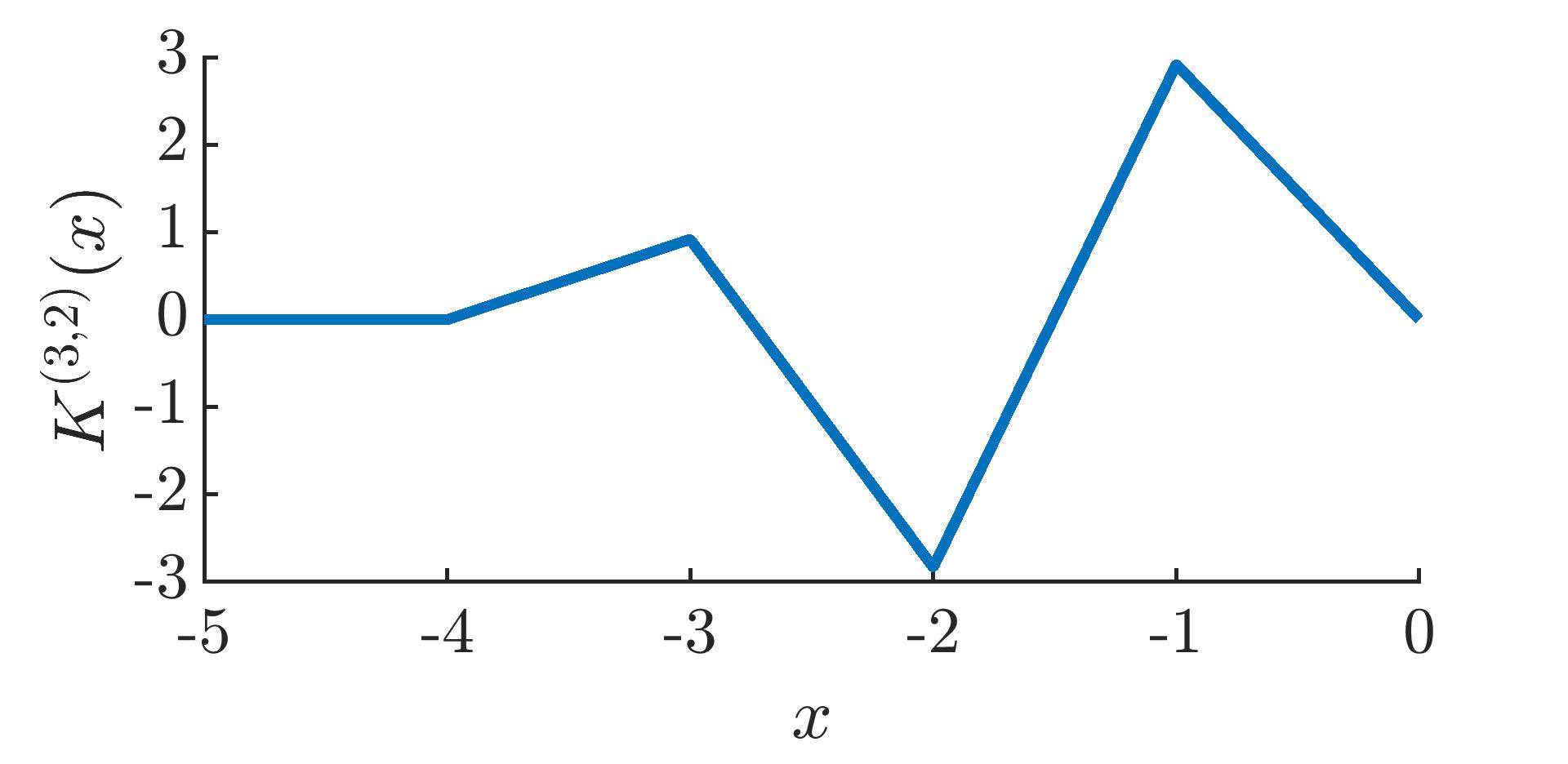}
  &     \includegraphics[width=0.33\linewidth]{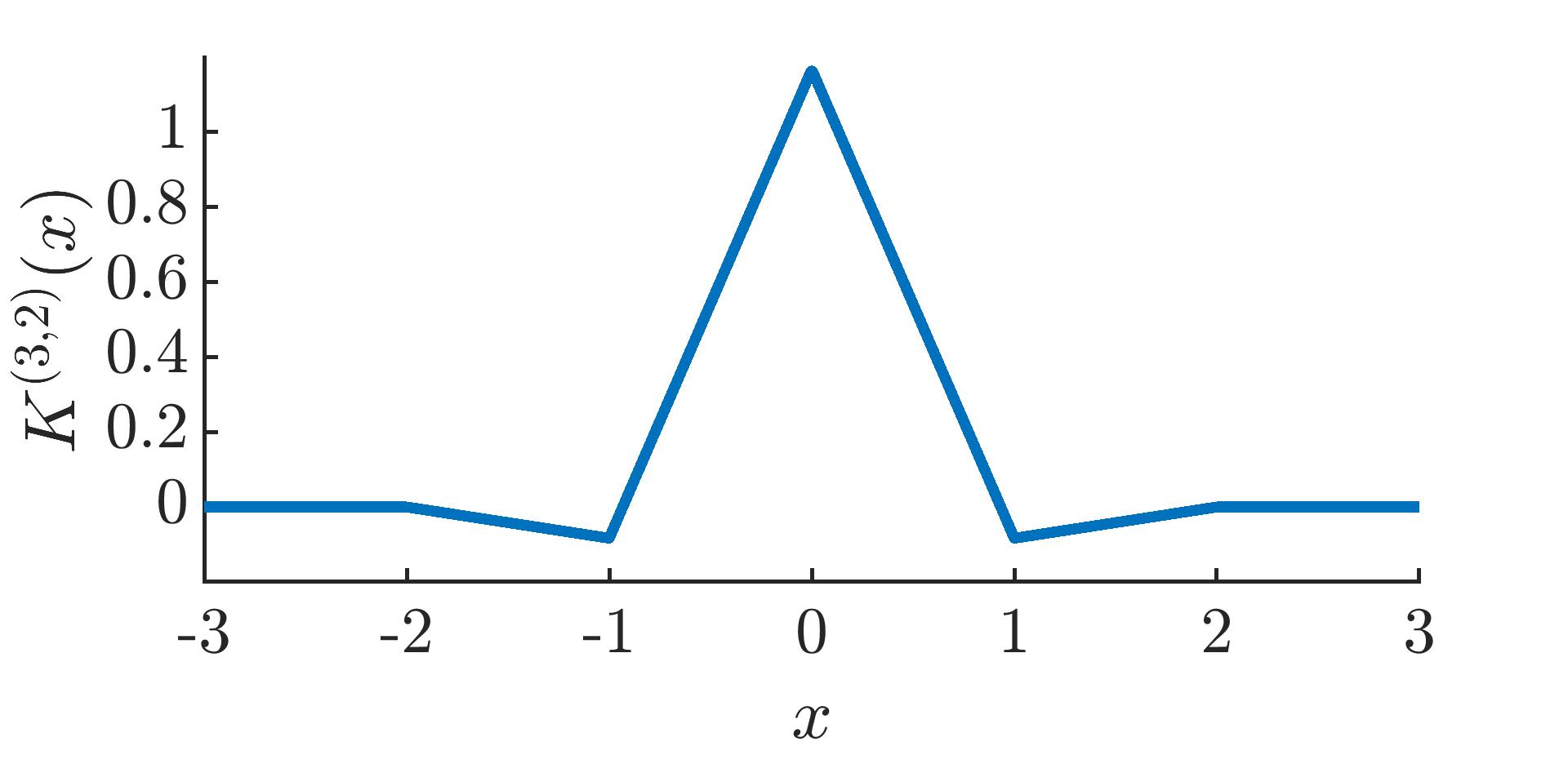}
     &     \includegraphics[width=0.33\linewidth]{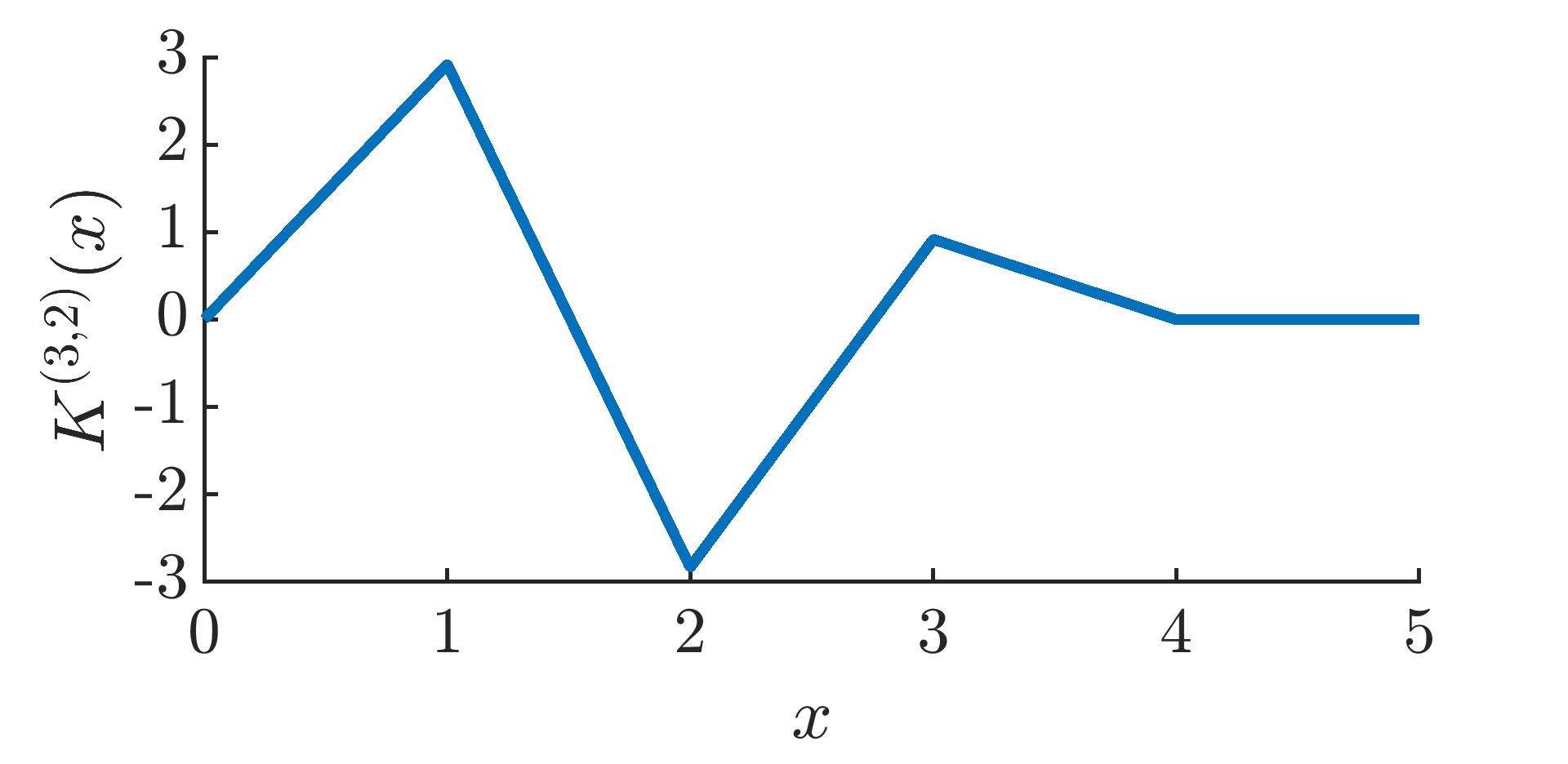}\\ \hline
    \end{tabular}
    \caption{Shifted Kernels for $K^{(3,2)}(x)$. For $\lambda=\pm2$, $x=0$ corresponds to the boundary. Note that owing to the reflection of $y$ in the kernel argument during the convolution, $\lambda=-2$ ($\lambda=2$) corresponds to the left-sided (right-sided)  kernel in the physical space.}
    \label{fig:shift_kernel}
\end{figure}

Given $\Omega=[a,b]$, consider the shifting function originally proposed in \cite{Slingerland}:
\begin{equation}
    \lambda(x^{\star})=\begin{cases}
\min\left\{0,-(r+\ell)/2+(x^{\star}-a)/H(x^{\star}) \right\}, & x^{\star}\in [a,(a+b)/2),\\
\max\{0,(r+\ell)/2+(x^{\star}-b)/H(x^{\star})\}, & x^{\star}\in [(a+b)/2),b],
\end{cases}
\end{equation}
where $H$ is the kernel scaling function. If  $\lambda(x^\star)$ {is plotted} over the domain $[-1,1]$, 
as depicted in Fig.~\ref{fig:shift},  it is observed that the support of the kernel is being shifted away from the boundaries when getting close to both boundaries, while in the interior  the kernel support remains  symmetric. To understand the effect of the shift function on the kernel support,  a plot {of} the fully shifted left ($\lambda=-2$), right ($\lambda=2$), and symmetric kernels 
are given in Fig.~\ref{fig:shift_kernel}. The shifting function provides a continuous transition between the shifted 
and symmetric kernels, though the limited smoothness of the shifting function restricts the smoothness of the filtered 
approximation at the interface of shifted and unshifted kernels. 

\subsection{Knot Matrix Examples}
\label{subsec:knotmatrix}
Typical examples of knot matrices corresponding to the case of  $r+1=3$ B-splines of order $\ell=2$ are given by
\[
\mathbf{T}_{\rm Left}=\begin{bmatrix}
-4 & -3 & -2\\
-3 & -2 & -1\\
-2 & -1 & 0
\end{bmatrix},\;\;\; \mathbf{T}_{\rm Sym}=\begin{bmatrix}
-2 & -1 & 0\\
-1 & 0 & 1\\
0 & 1 & 2
\end{bmatrix},\;\;\; \mathbf{T}_{\rm Right}=\begin{bmatrix}
 0 & 1 & 2\\
1 & 2 & 3\\
2 & 3 & 4
\end{bmatrix},
\]
where $\mathbf{T}_{\rm Left}$ and $\mathbf{T}_{\rm Right}$ denote the biased filters used for filtering at the left and right boundaries of the domain, respectively. The symmetric kernel $T_{\rm Sym}$ is used in the interior of the domain.

\subsection{Kernel coefficients}
\label{subsec:coeff}

To determine the coefficients, $c_{\gamma},$ the post-processing filter must satisfy consistency and moment conditions. These are equivalent to polynomial reproduction:
\begin{center}
\bgroup
\def\arraystretch{1.5}
\label{consistency}
\begin{tabular}{l c l}
{\sf Consistency + $r$ Moments} & & {\sf Polynomial Reproduction}\\
$\int_{\mathbb{R}}\, K_{\mathbf{T}}(x)\, dx = 1$ & $\Rightarrow$ & $\int_{\mathbb{R}}\, K_{\mathbf{T}}(x-y)\, dy = 1$ \\
$\int_{\mathbb{R}}\, K_{\mathbf{T}}(x)x^k\, dx = 0, \quad 1\leq k \leq r$ & $\Rightarrow$ & $\int_{\mathbb{R}}\, K_{\mathbf{T}}(x-y)y^k\, dy = x^k,\, k\in \mathbb{Z}, \quad 1\leq k \leq r$
\end{tabular}
\egroup
\end{center}

\noindent where $r$ is the number of moments. This aids in ensuring that the accuracy of the underlying data is not destroyed. For non-symmetric kernels the knot matrices can change with each filtering point $x^{\star}$, and so the kernel coefficients need to be recomputed with every convolution evaluation, which is not the case in the interior. The polynomial reproduction requirement
\[K_{\mathbf{T}}\star x^p=x^p,\;\;\;p=0,\hdots, r,\]
can be expressed as
\[
\sum_{\gamma=0}^rc_{\gamma}(x^{\star})\int_{\mathbb{R}}B^{\ell}_{\mathbf{T}_{\gamma}}(y)(x-y)^p\;dy=x^p,\;\;\;p=0,\hdots,r.
\]
  Choosing $x=0$ results into a linear system
    \begin{align}
   \label{coeff_system}\begin{bmatrix}
    b_{0,0}&\hdots & b_{0,r}\\
    \vdots& &\vdots\\
    b_{r,0}&\hdots & b_{r,r}
    \end{bmatrix}\begin{bmatrix}
    c_0\\ c_1\\ \vdots \\ c_{r}
    \end{bmatrix}=\begin{bmatrix}
    1\\0\\ \vdots \\ 0
    \end{bmatrix} \end{align}
    where
    \[b_{p,\gamma}=\int_{\mathbb{R}} B^{\ell}_{\mathbf{T}_{\gamma}}(y) (-y)^p\;dy\;\;\text{for}\;\;p,\gamma=0,\hdots,r .\]
Solving the linear system~\eqref{coeff_system} yields the kernel coefficients. More explicit formulas for the  coefficient system \eqref{coeff_system} in the uniform symmetric knot matrix case are provided in \cite{Mirzargar2016}, while the general knot matrix case is considered in \cite{peters2015,peters2016}.

\subsection{Generalized Spline Kernels}
An issue arising from position-dependent filters is that errors increase in boundary regions for FEM data and in the case of PIC denoising, 
 resulting into the failure to preserve boundary layers. In the past these issues have been addressed by using larger numbers of B-splines or convex combinations of SIAC kernels \cite{Slingerland,Ji2014}, but as detailed in \cite{Ryan2015}, adding a single generalized spline can decrease errors at the boundaries without increasing the kernel support. For example, near the left boundary, write
\[
K_{H \, \mathbf{T}}(x)=\sum_{\gamma=0}^{r}c_{\gamma}(x) B^{\ell}_{H \, \mathbf{T}_{\gamma}}(x)+c_{r+1}B^{\ell}_{H\, \mathbf{T}_{r+1}}(x).
\]
Here the coefficients are still determined  by requiring polynomial reproduction. The new knot matrices are defined based off the sign of $\lambda$ by
\[
\lambda<0 \rightarrow \begin{bmatrix}
\mathbf{T}\\
\mathbf{s}_{L}
\end{bmatrix},\;\; \lambda=0 \rightarrow \mathbf{T},\;\; \lambda>0 \rightarrow \begin{bmatrix}
\mathbf{s}_R\\
\mathbf{T}
\end{bmatrix}.
\]
Here the $1\times (\ell+1)$ knot sequences for the left and right generalized splines are given by
\[
\mathbf{s}_L=\Big[\lambda+\frac{1}{2}(r+\ell)-1, \lambda+\frac{1}{2}(r+\ell), \hdots, \lambda+\frac{1}{2}(r+\ell)\Big],
\]
and
\[
\mathbf{s}_R=\Big[\lambda-\frac{1}{2}(r+\ell), \hdots, \lambda-\frac{1}{2}(r+\ell), \lambda-\frac{1}{2}(r+\ell)+1\Big].
\]

\hfill

\begin{figure}
    \centering

\begin{tabular}{|c|c|}\hline
\multicolumn{1}{|c|}{$\lambda=-2$}&\multicolumn{1}{|c|}{$\lambda=2$}\\ \hline
       \includegraphics[width=0.33\linewidth]{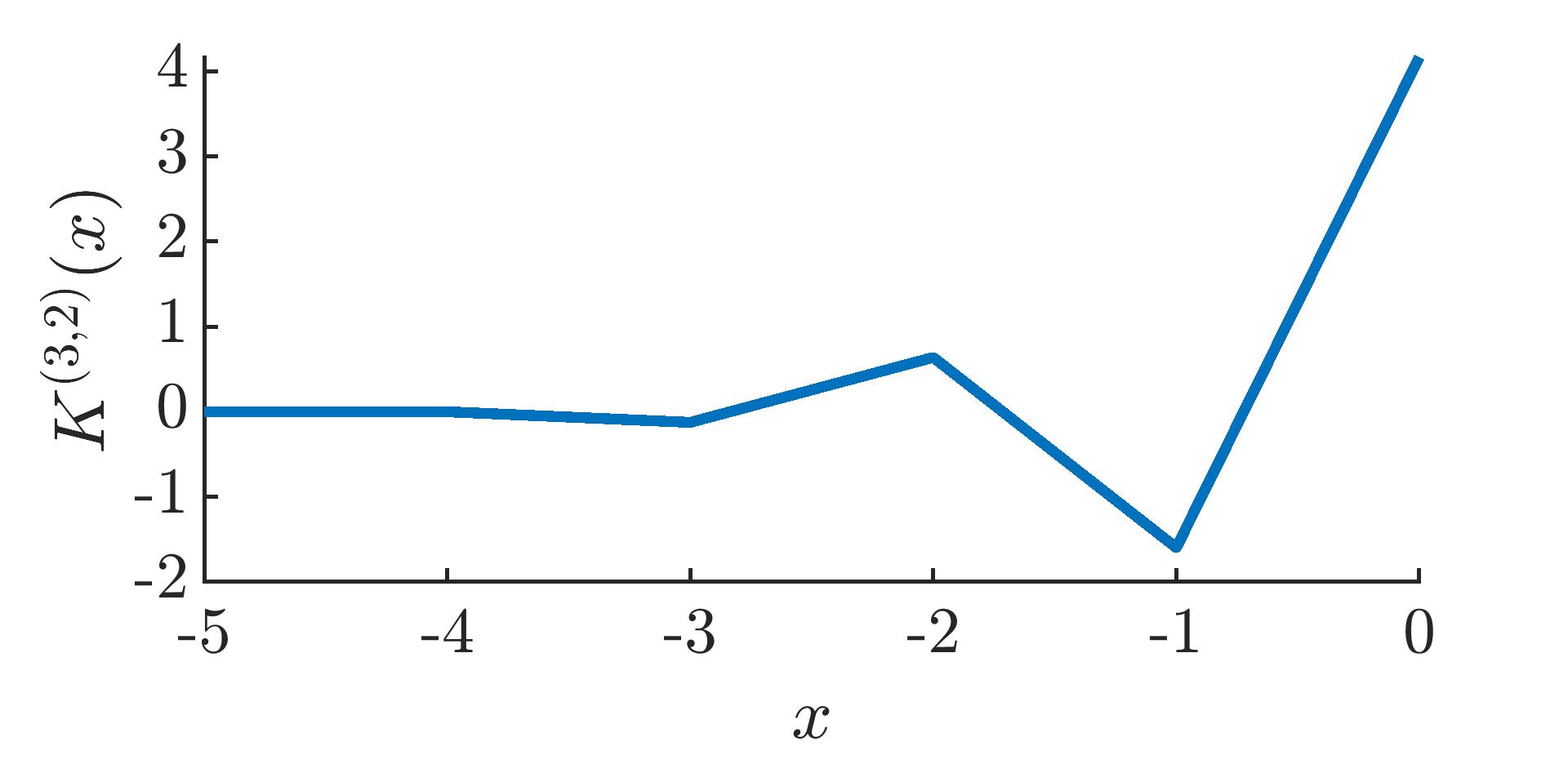}
     &     \includegraphics[width=0.33\linewidth]{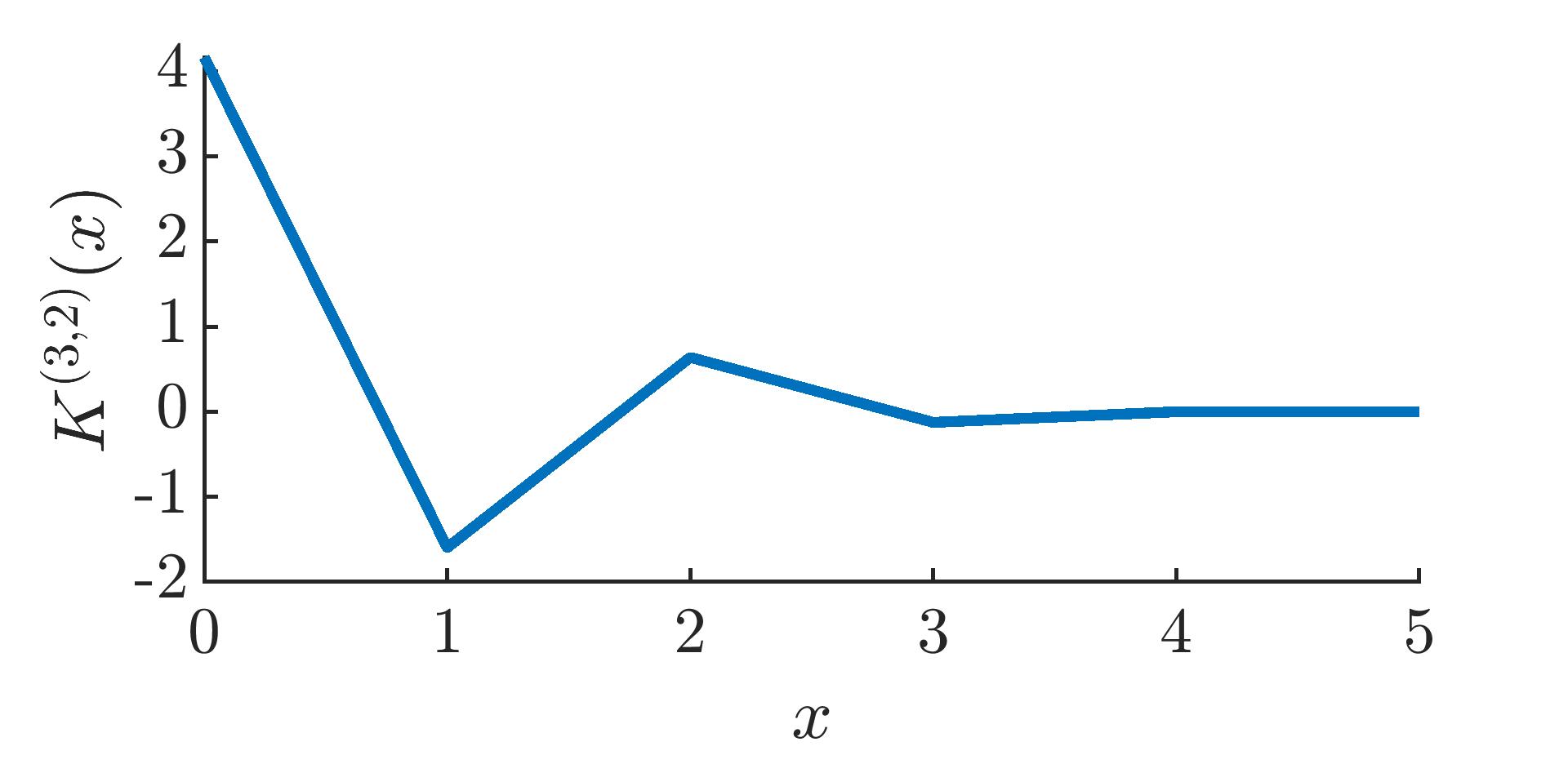}\\ \hline
    \end{tabular}

    \caption{Shifted Kernels for $K^{(3,2)}(x)$ including a generalized spline. Here $x=0$ corresponds to a boundary.}
    \label{fig:shift_kernel_gen}
\end{figure}

For instance, consider the case of $r=2$ and $l=2$. The new left and right knot matrices for $K^{(3,2)}_H$ are given by
\[
\mathbf{T}_{Left}=\begin{bmatrix}
-4 & -3 & -2\\
-3 & -2 & -1\\
-2 & -1 & 0\\
-1 & 0 & 0
\end{bmatrix},\;\;\;  \mathbf{T}_{Right}=\begin{bmatrix}
 0&0&1\\
 0 & 1 & 2\\
1 & 2 & 3\\
2 & 3 & 4
\end{bmatrix}.
\]
Plots of the these boundary kernels are given in Fig.~\ref{fig:shift_kernel_gen}.

\subsection{Adaptive kernel scaling}
\begin{figure}
    \centering
    \includegraphics[width=0.75\linewidth]{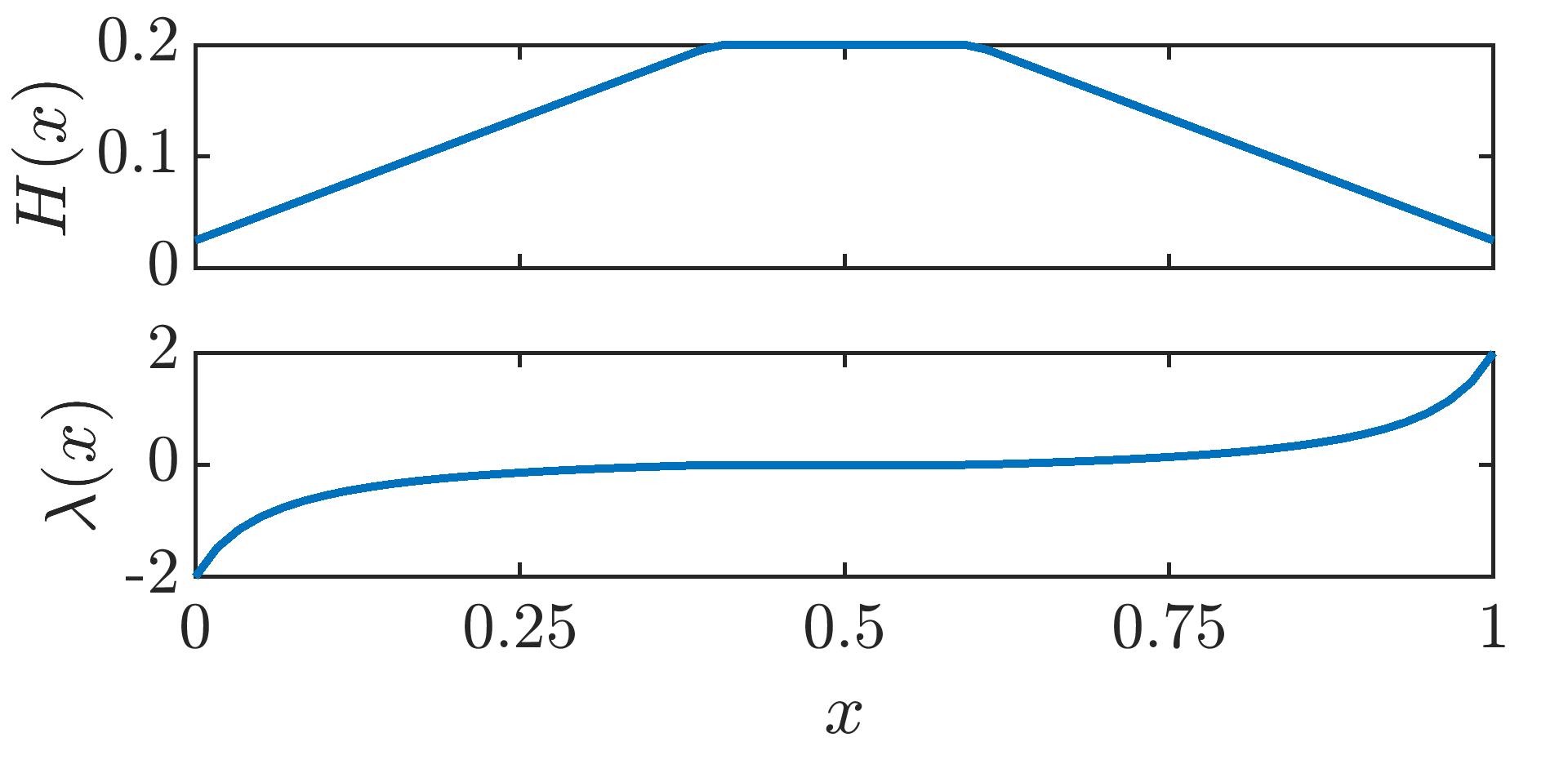}
    \caption{An adaptive kernel scaling example for $K^{(3,2)}(x)$ over the domain $\Omega=[0,1]$ with $H_{\rm int}=0.2$, $h_{\rm grid}=0.1$.}
    \label{fig:adaptive}
\end{figure}

Although the addition of generalized splines improves the performance of SIAC denoising for moderate kernel scalings, as will be observed in the numerical results, 
some care is needed to preserve short-duration sharp gradients near boundaries when the kernel scaling becomes large. To that end, we propose a remedy via an adaptive kernel scaling   given by
\begin{equation} \label{eqn:adapt_scale}
H_{\rm adaptive}(x)=\begin{cases}
H_{\rm int},\;& x-H_{\rm int}\frac{r+\ell}{2}\geq a \;\text{ and }\;x+H_{\rm int}\frac{r+\ell}{2}\leq b,\\
\frac{h_{\rm grid}}{r+\ell}+2(x-a)\frac{H_{\rm int}-h_{\rm grid}/(r+\ell)}{H_{\rm int}(r+\ell)},\;& x-H_{\rm int}\frac{r+\ell}{2}\leq a,\\
\frac{h_{\rm grid}}{r+\ell}+2(x-b)\frac{-H_{\rm int}+h_{\rm grid}/(r+\ell)}{H_{\rm int}(r+\ell)},\;& x+H_{\rm int}\frac{r+\ell}{2}\geq b,
\end{cases}
\end{equation}

\noindent where $H_{\rm int}$ is a large scale desired for the interior, and $h_{\rm grid}$ is the size of the mesh elements. 
$H_{\rm int}$ is typically related to the physical length scale in the problem. 
This choice of adaptive scaling provides the advantage of preserving the boundary values of the underlying initialization by reducing the kernel support to a single element at the boundary locations. This enables the polynomial-reproduction property of the kernel to be maintained. 
Fig.~\ref{fig:adaptive} shows a depiction of this adaptive scaling and its smoothing effect on the transition between {the symmetric kernel and the} 
position-dependent kernel with a shifting parameter. This adaptive scaling has the additional effect of creating a smoother shifting function $\lambda$ and thereby resulting in a smoother transition between different knot matrices.

\subsection{Convolution evaluation}
The consistent evaluation of convolution in 1D requires keeping track of all the breaks in continuity of the approximation and breaks in differentiability of the kernel within the scaled kernel support. By only integrating between these breaks in continuity of the data or the kernel,  the convolution 
{can be performed} exactly using Gaussian quadrature. A detailed discussion of the implementation can be found in Docampo-S\'anchez et al. \cite{whatWorks}.

\subsection{Fourier representation}
\label{sec:fourier}
As our primary focus is to denoise PIC data,  it is critical to examine the Fourier transform of the convolution kernel to quantify the damping effect applied to each frequency.
This examination will guide the choice of the kernels when denoising moments. 
Unfortunately, analytical computation of the kernel's Fourier response is only possible 
 in the periodic case for the symmetric kernel with a constant kernel scaling. 
 In such a case, the Fourier representation of the SIAC kernel composed of $(r+1)$ B-splines of order $\ell$ and constant scaling $H$ is given by
 \begin{equation}
 \label{eq:KFT}
\hat{K}^{(r+1,\ell)}_H(k)=\left(\frac{\sin(\frac{kH}{2})}{\frac{kH}{2}}\right)^{\ell}\left(c_{{\lceil\frac{r}{2}\rceil+1}}+2\sum_{\gamma=1}^{{\lceil\frac{r}{2}\rceil}}\, c_\gamma\cos(\gamma kH)\right).
 \end{equation}
 The first term in the product is controlled by the spline order, while the second is controlled by the number of B-splines, 
 {with the value of the coefficients, $c_\gamma$, being controlled by both the number of B-splines and the B-spline order}.  
 
 \hfill
 
\begin{figure}[tp!]
      \centering
 \begin{tabular}{c c c}
     \includegraphics[width=0.45\linewidth]{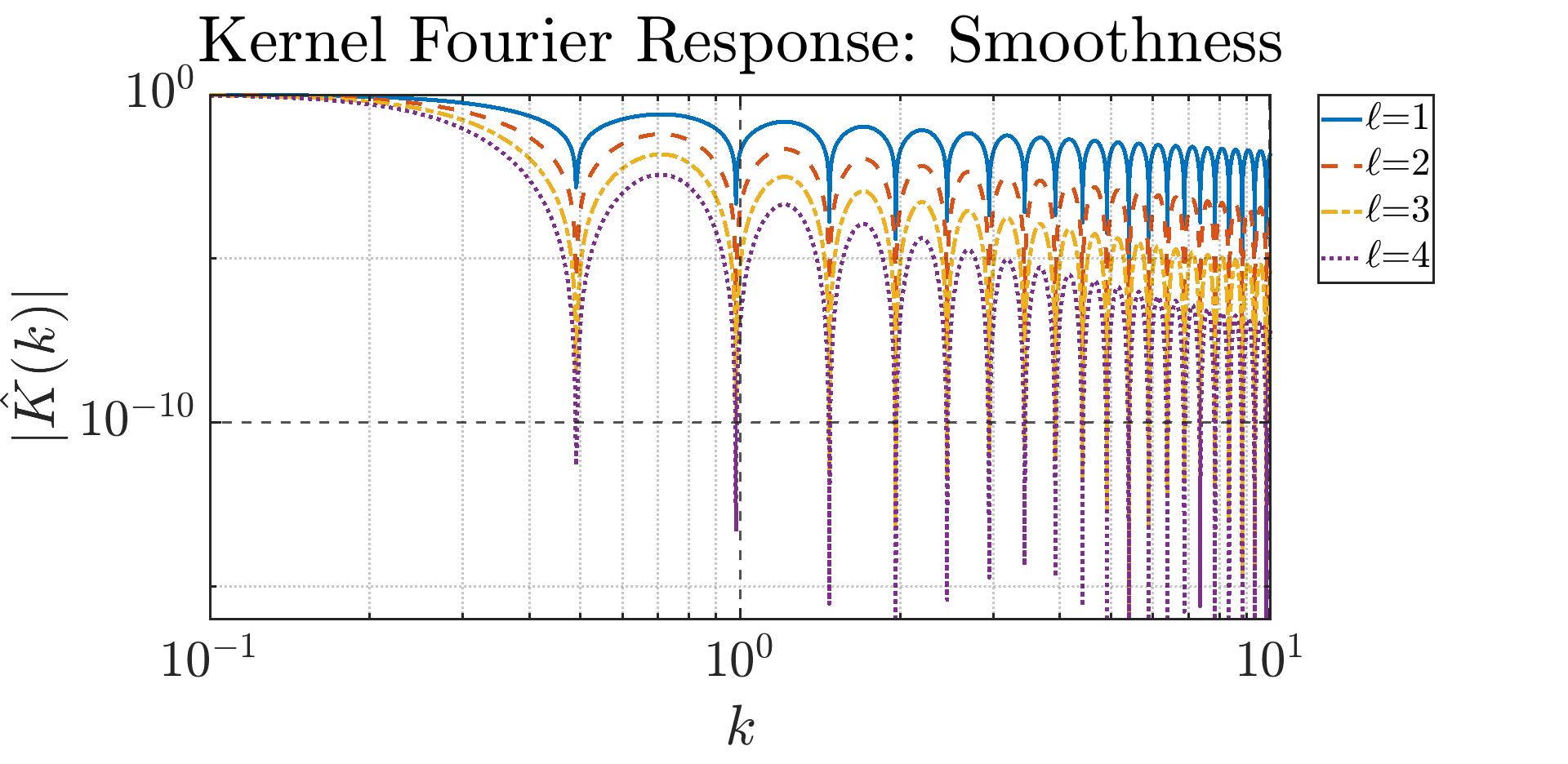}
& &

          \includegraphics[width=0.45\linewidth]{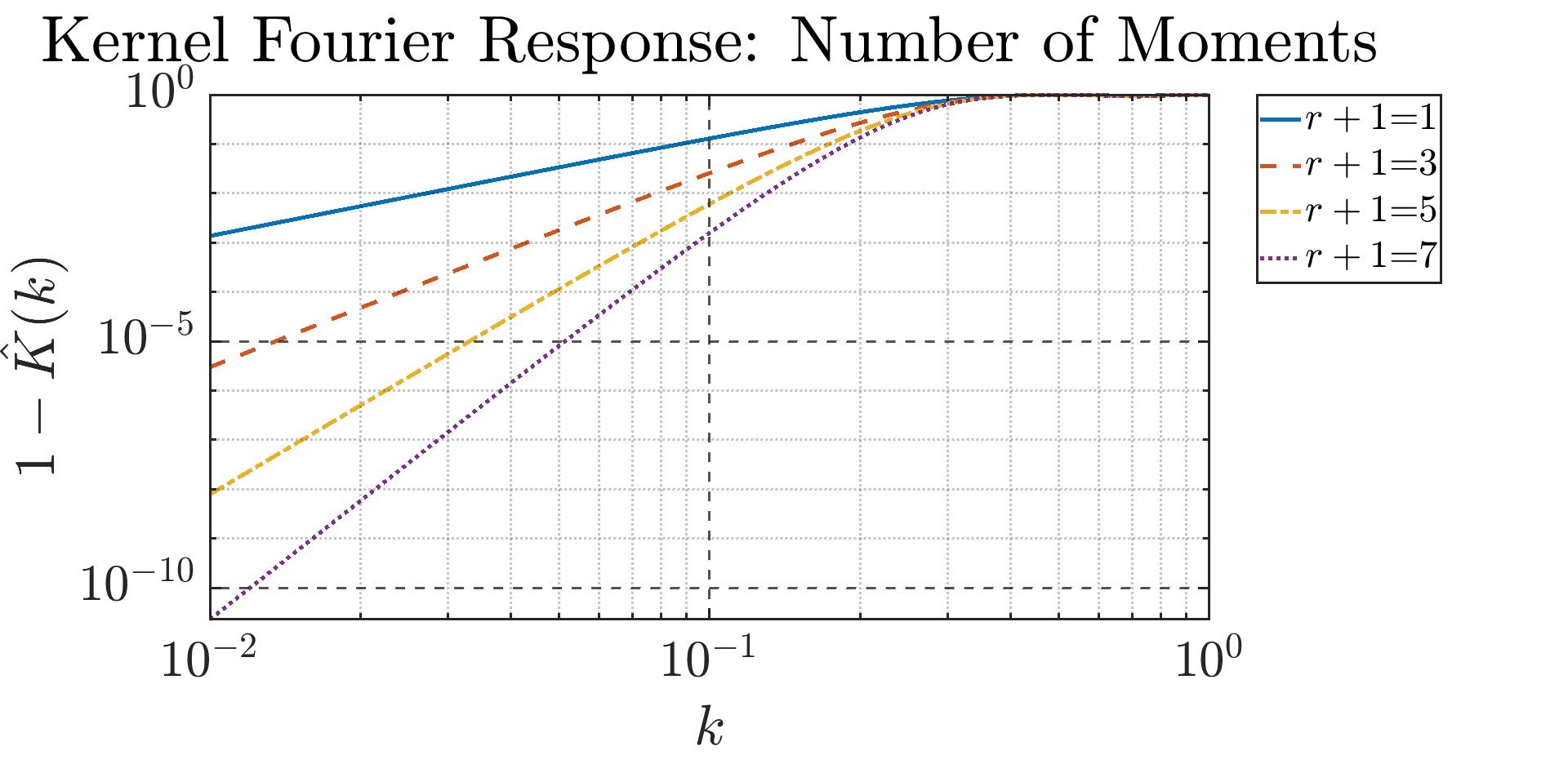}\\
          (a) && (b)
      \end{tabular}
\centering
     \includegraphics[width=0.45\linewidth]{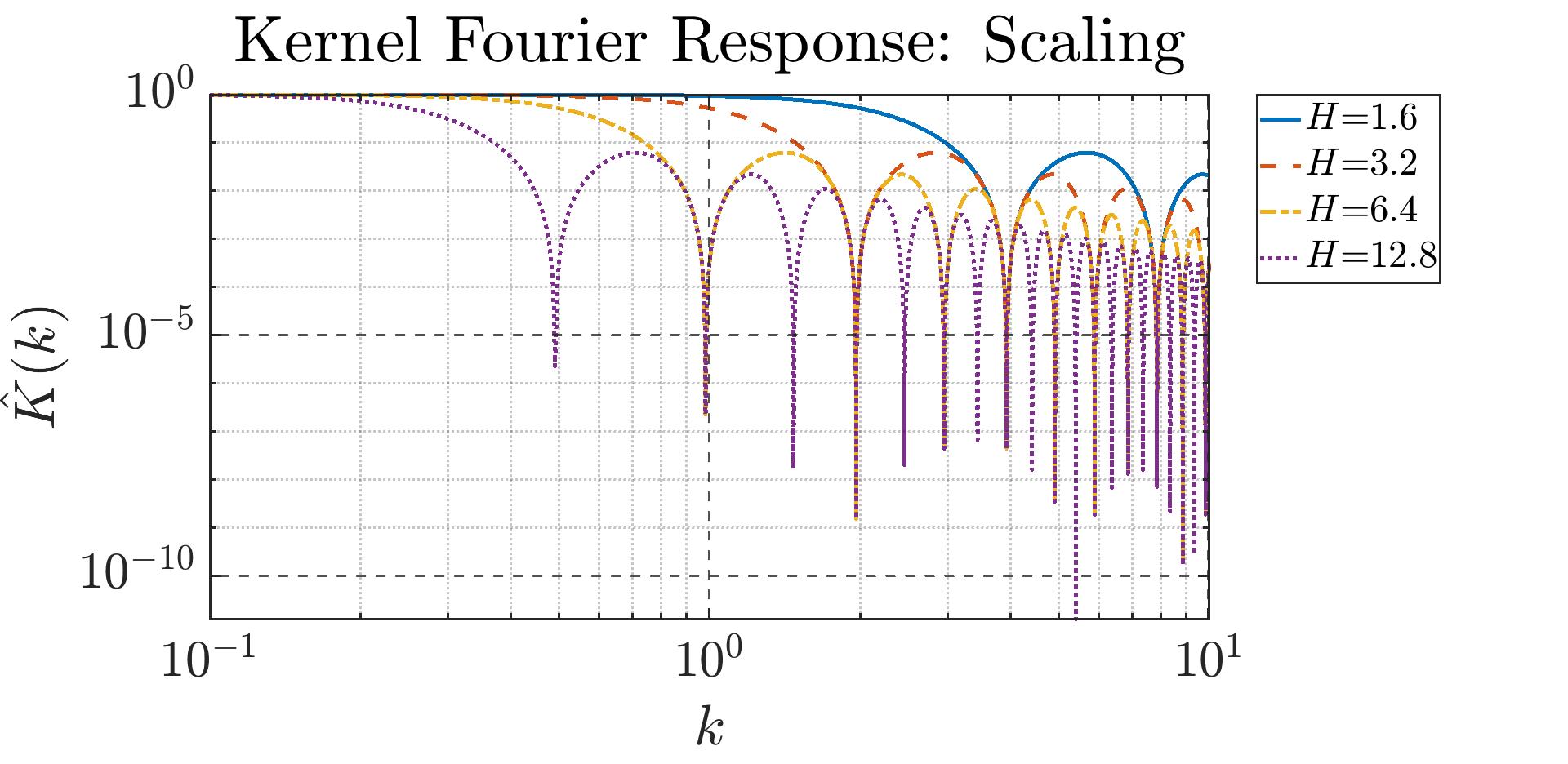}\\
     (c)

    \caption{Variation of the magnitude of the analytic Fourier representation of the SIAC filter, $K^{(r+1,\ell)}_H$, with respect to (a) spline order ($\ell$), with higher orders corresponding to smoother and more diffused filtered data, (b) moment conditions ($r$), with more moment conditions allowing the capturing of more information about a point, and (c) kernel scaling ($H$), with a larger scaling resulting in more information being used in the filtering process. Note that the kernel is symmetric with respect to $k$ and that, unless indicated otherwise in the legends, $H=12.8$, $r+1=3$, and $\ell=2$.}
    \label{fig:Fourier} 
 \end{figure}

 In Fig.~\ref{fig:Fourier} the Fourier representation is displayed for varying kernel parameters. We observe that increasing the spline order, which corresponds to an increase in smoothness, causes the filter to sample smaller magnitude wave numbers more heavily and reduce oscillations. This makes sense as the Fourier transform of smoother functions decays faster and increasing the spline order results in a smoother kernel. Increasing the number of splines, which corresponds to increasing the number of moments, does not change the shape drastically, but slightly alters how the kernel samples different modes. Specifically, it flattens the Fourier response at $k=0$ by causing $\frac{d^n}{dk^n}\hat{K}^{(r+1,\ell)}(0)=0$ for $n=1,\hdots,r$. This follows from the polynomial reproduction condition on the kernel.  Lastly, changing the scaling factor only causes a dilation in $k$ of kernel response and utilizes more (or fewer) elements in the filtering process. An understanding of the frequency effects of position-dependent and adaptively scaled filters remains an important area of further study. In particular, the standard methodology of computing the Fourier transform of the SIAC kernel is inappropriate owing to the spatial variability of the kernel function employed over the computational domain. One possible approach is the consideration of the spectrum of the discrete filtering operator, but that is beyond the scope of this work.
 
\section{Numerical Results\label{Numerical Results}}

In this section the applications of the SIAC filter to
denoise PIC data {are considered.} Here,  both a periodic boundary condition and a boundary condition that
introduces sharp gradients will be considered. For the former, an investigation of the
damping of the noise using the uniformly scaled symmetric filter is done, where it is expected that the filter
 should preserve the physical oscillations. For the latter case, several aspects are studied, including
comparisons of the abilities of position-dependent filters with or
without generalized splines as well as with or without adaptive kernel scalings. 
The focus is if the boundary layers are resolved in the physical space and how that affects the Fourier space. 
It must be highlighted again that maintaining
the boundary behavior of the data can be critical to obtaining
accurate physics, e.g., for the Bohm criterion.


\newpage

Given PIC moment data, there are three steps in the denoising procedure:
\begin{enumerate}[leftmargin=2cm, label=Step \arabic*:]
    \item Construct a nodal piecewise interpolant as described in Sec.~\ref{Intialization} for the given pointwise discrete data {that is generated} from {a} PIC simulation. 
    \item Select the kernel parameters {($r,\, \ell,$ and $H$)}, and filter the initialization obtained in Step 1 with the corresponding SIAC kernel. The choice of parameters depends on the type of data being filtered (characteristic length scales, desired damping effects, etc.).
    \item Assess the filtered data from the point of view of physics
\end{enumerate}
It was determined that there is an insensitivity in the initialization step to low polynomial degrees and therefore only results for a piecewise constant initialization are included in what follows.

\paragraph{Assessment of Fourier Effects}
For the purposes of displaying the Fourier effects of the filtering procedure on the underlying data,  plots of the single-sided amplitude spectrum {are given}. 
If the data is non-periodic,  multiplication of the data vector by {a} Hanning window function {is done} in the first step to minimize any spectral leakage.

\subsection{Periodic boundary conditions}
\label{ex:periodic}
\begin{figure}[tp!]
      \centering
\includegraphics[width=0.9\linewidth]{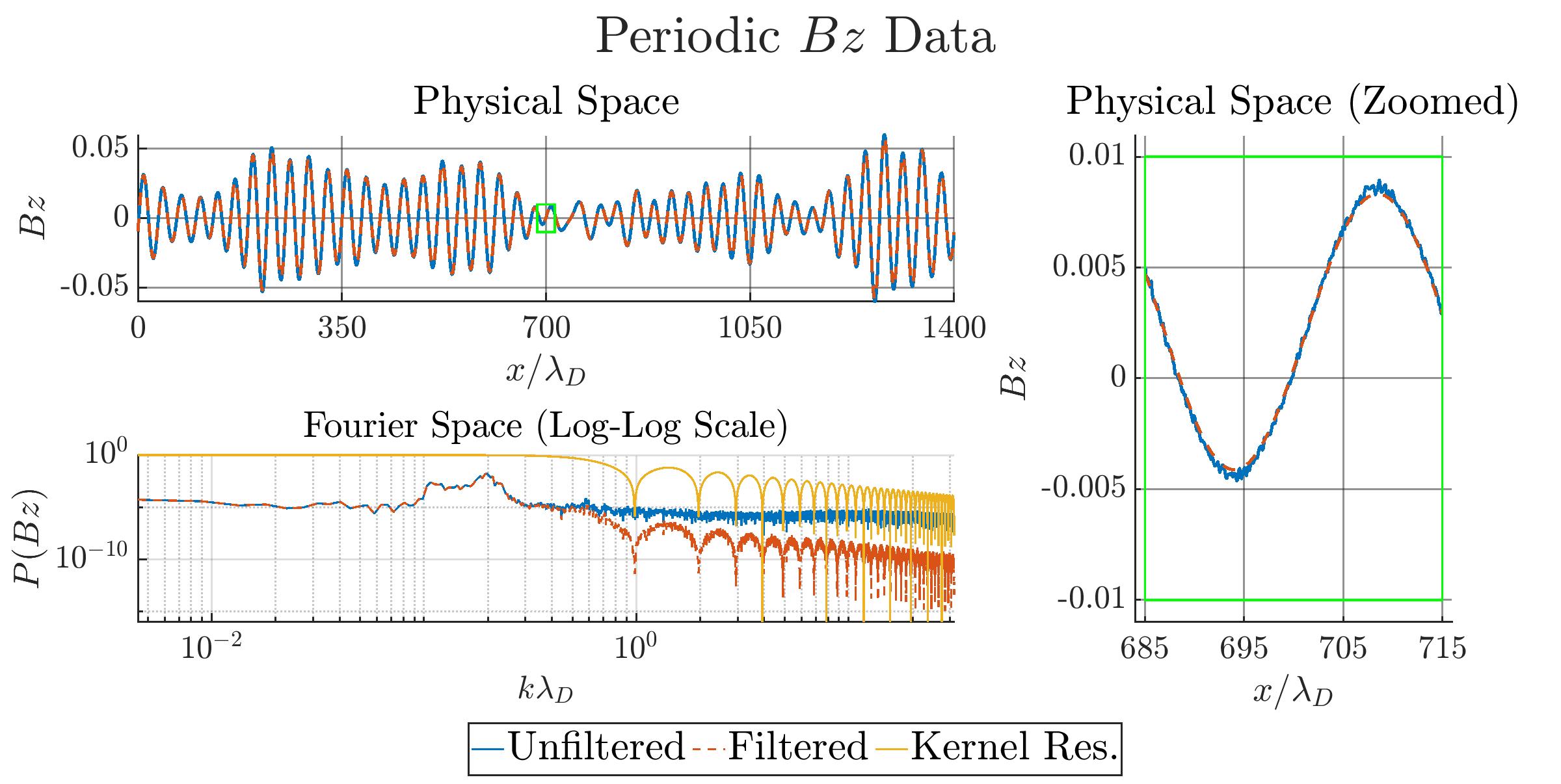}
          \caption{Application of SIAC filtering to periodic magnetic field for the whistler instability driven by the trapped electrons: (upper-left) filtered vs. unfiltered data over the whole domain. (right) zoomed-in physical space plot displaying spatial noise damping.  (lower-left) single-sided amplitude spectrum of filtered and unfiltered data with analytic SIAC kernel response superimposed. }
     \label{fig:magnetic_per}
      \end{figure}

 The first test {is for} the possible effectiveness of the SIAC technique on
denoising the PIC data. This requires the removal of high
frequency/wavenumber (with wavelengths at cell size) noises and the
preservation of physical oscillations. For this purpose, 
PIC simulations with periodic boundary conditions {are considered}. Such boundary
conditions can be widely found in PIC
simulations~\cite{birdsall2018plasma,hockney2021computer} such as in the
study of plasma
instabilities~\cite{pritchett2000particle,verboncoeur2005particle}. Here
 the data from the VPIC code~\cite{VPIC} that simulates the
whistler instability in magnetized plasma induced by the trapped
electrons {is used}~\cite{guo2012ambipolar,zhang2023collisional}. For such
periodic data, a symmetric SIAC filter with a constant
kernel scaling over the whole domain {can be used}. 
Fig.~\ref{fig:magnetic_per} shows  the results of applying the
filtering procedure on periodic magnetic field data where the kernel
is composed of 3 B-splines of order 2 ($r+1=3$, $\ell=2$). Here the
domain is $\Omega=[0,1400]$ with {a} uniform grid spacing {of} $h_{grid}=\Delta x =0.1$, in the
unit of Debye length $\lambda_D$.  
A large $H=64 h_{grid}$ {is chosen} so
that the truncated wavenumber is $k_c
\tilde{<}\lambda_D^{-1}$.Choosing a small filter scaling, $H$,  would result in less damping of the higher order modes as illustrated in Equation \eqref{eq:KFT} and Fig.~\ref{fig:Fourier}.


We first note that the Fourier modes match very well in the low frequency domain,
resolving all the physical structures we are interested in. 
It is even more remarkable that 
the theoretical model of the
kernel response matches very well {with} the observed damped frequencies in the Fourier space (the lower-left plot).
As the analysis in Sec.~\ref{sec:fourier} suggested, the damped frequencies are primarily controlled by the the kernel scaling $H$. 
Therefore, by choosing the kernel scaling
appropriately,  it is possible to target and damp noises with high
wavenumbers while a good preservation of the real whistler waves peaked at
$k\sim 0.1\lambda_D^{-1}$ {is maintained}, as seen from the Fourier spectrum.

 This result suggests that for periodic data with limited
sharp features SIAC filtering {can be applied} in an informed manner.
The targeting physical frequency such as the whistler wave case is typically known {\it a prior} and 
can be readily used to guide the choice of $H$ in practice. 



 

      
\subsection{Data with sharp gradients near the boundary}
\label{ex:gradients}
\begin{figure}[tp!]
\begin{tabular}{|c|c|}\hline
  \multicolumn{1}{|c|}{$H=3.2$}   & \multicolumn{1}{|c|}{$H=12.8$} \\ \hline
    \multicolumn{1}{|c|}{Without generalized splines, constant scaling}   & \multicolumn{1}{|c|}{With generalized splines, constant scaling} \\ \hline
     \includegraphics[width=0.45\linewidth]{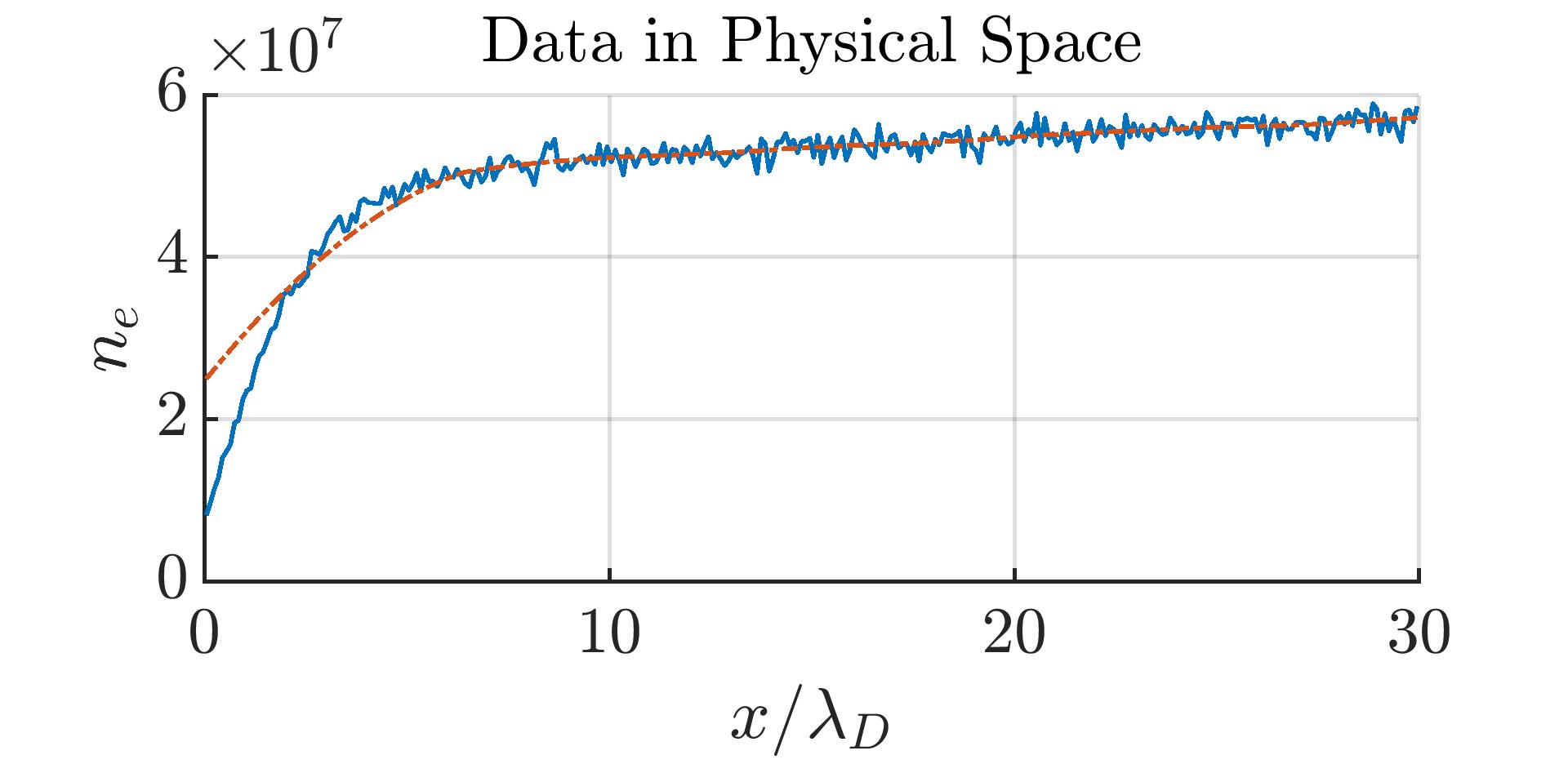} &  \includegraphics[width=0.45\linewidth]{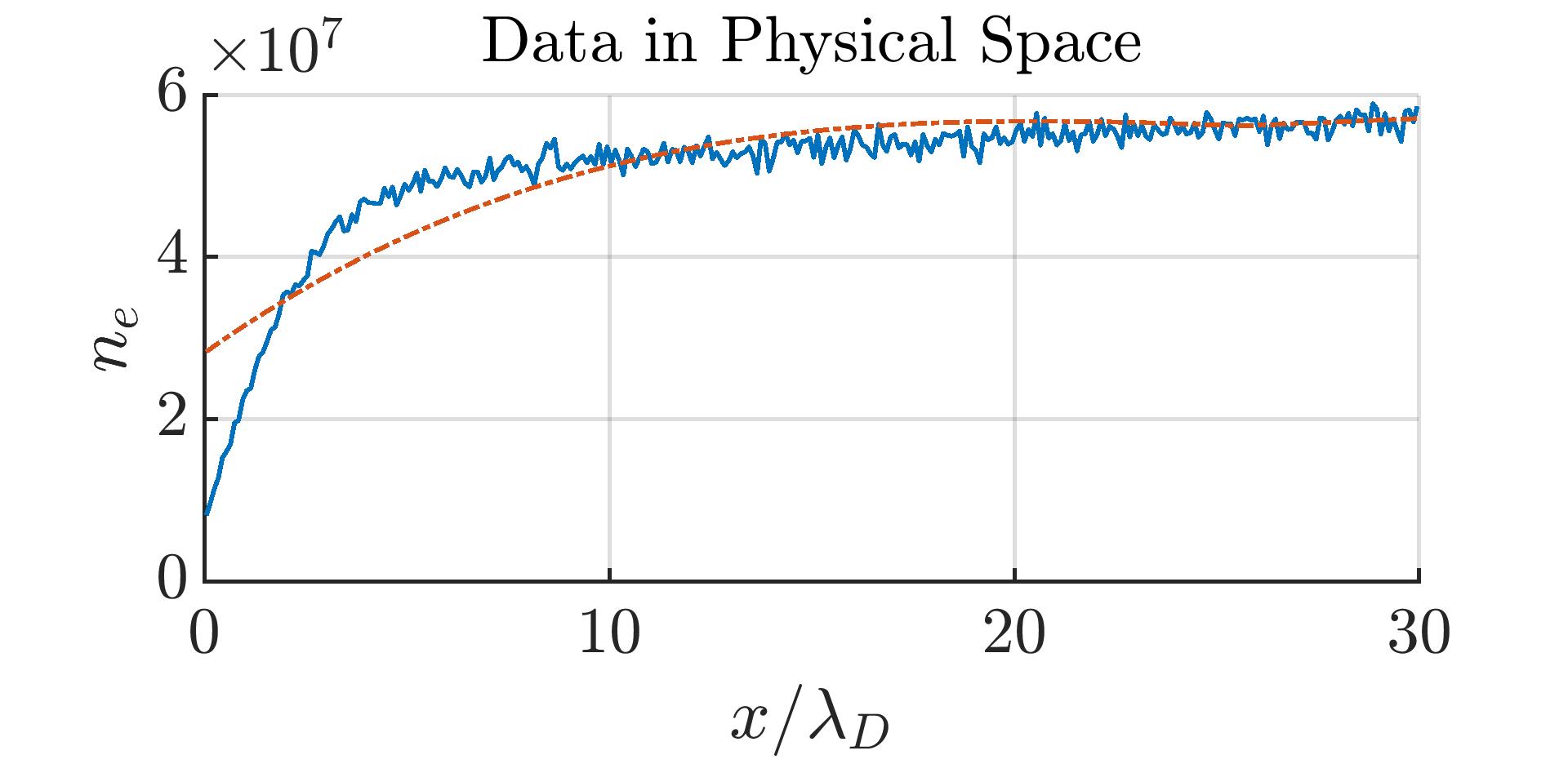}\\ \hline
         \multicolumn{1}{|c|}{With generalized splines, constant scaling}   & \multicolumn{1}{|c|}{With generalized splines and adaptive scaling} \\ \hline
      \includegraphics[width=0.45\linewidth]{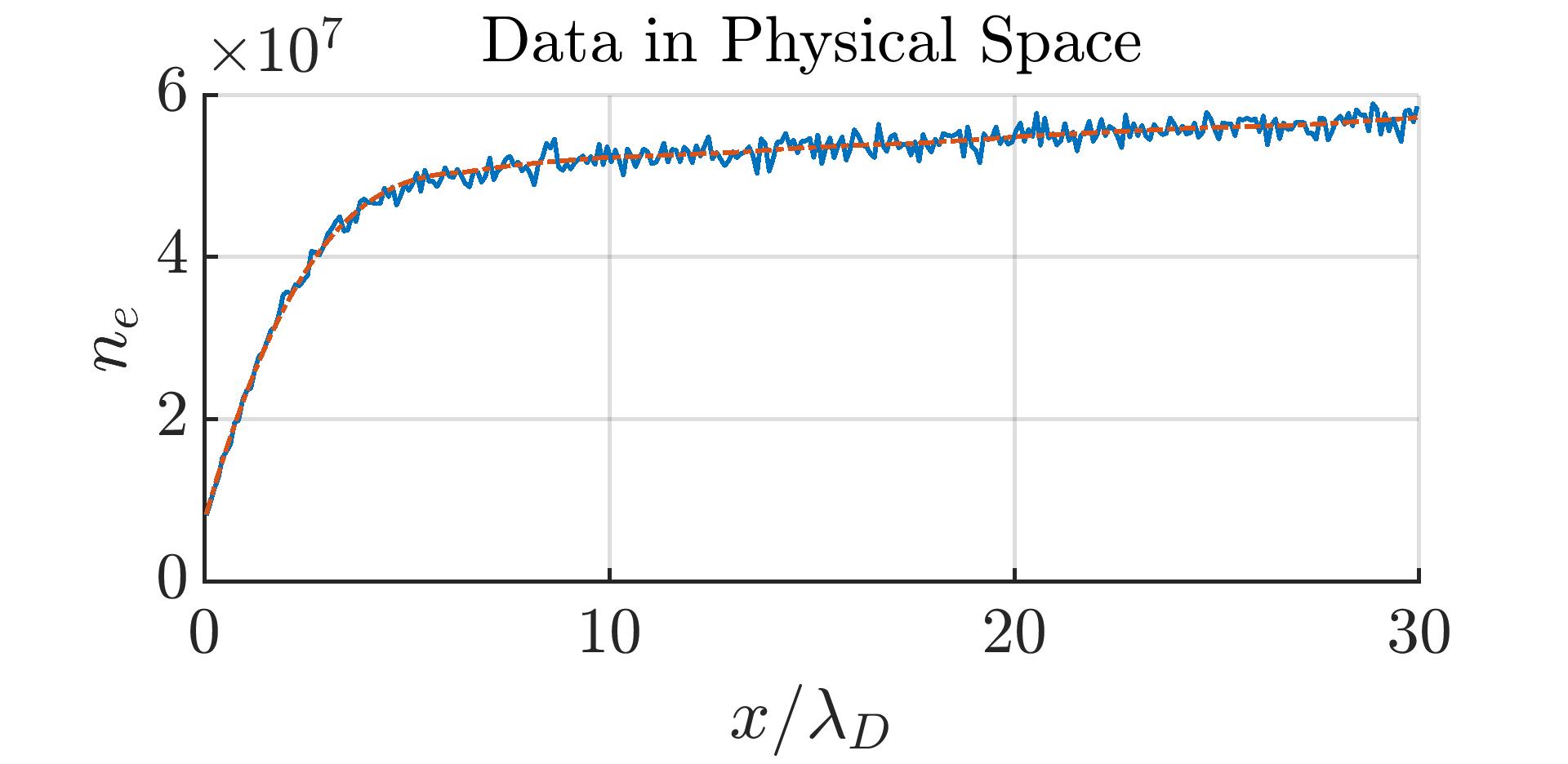} &  \includegraphics[width=0.45\linewidth]{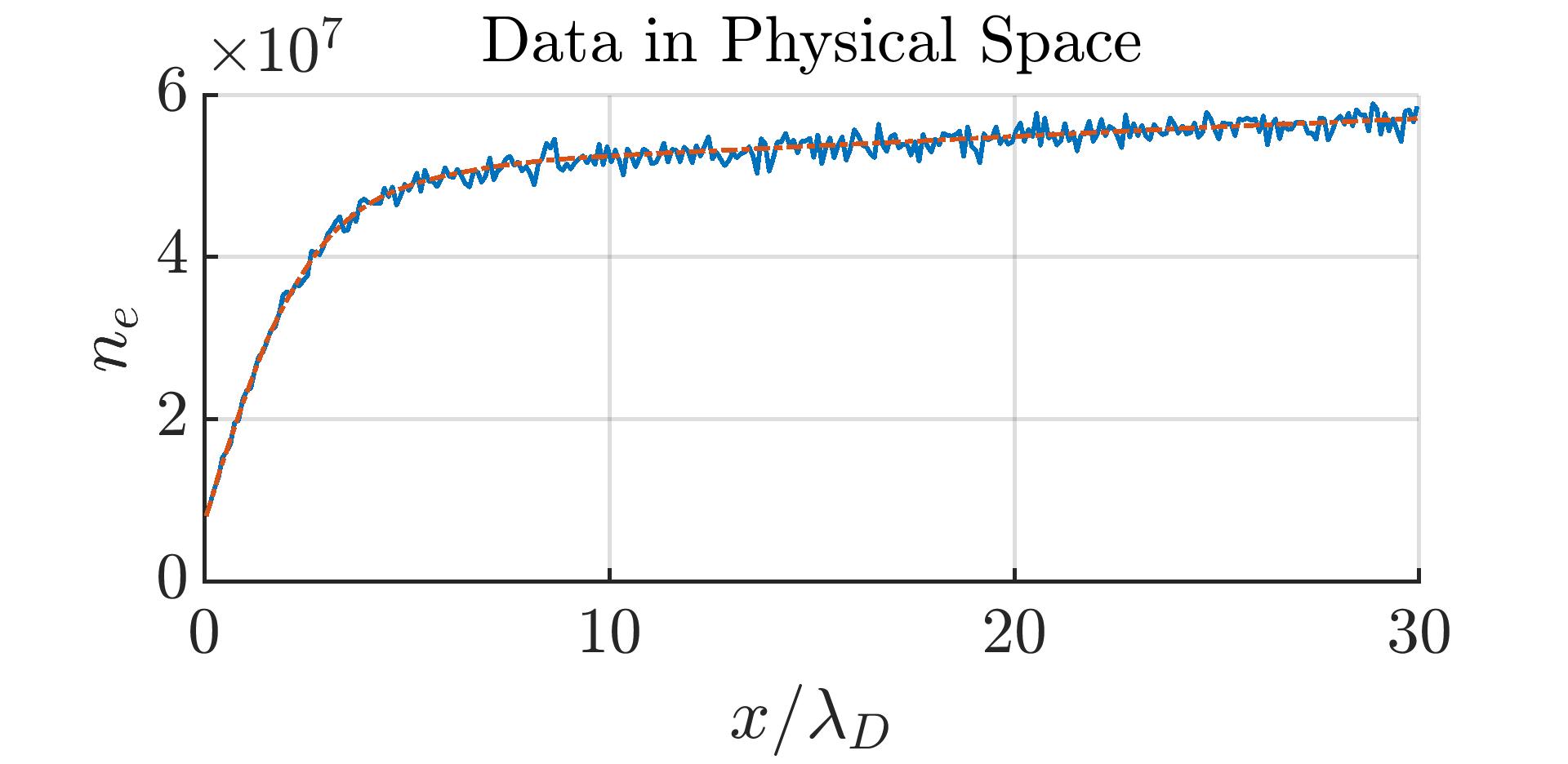}\\ \hline
\end{tabular}
    \centering
    \caption{\label{fig:bnd_behavior_ne} Effect of generalized splines and adaptive kernel scaling on boundary behavior preservation. The underlying data is electron density data $n_e$ with $\Omega=[0,800]$ and $h_{grid}=\Delta x=0.1$. Here the solid blue lines denote the unfiltered data and the dashed red lines the filtered data.}
\end{figure}

\begin{figure}[tp!]
\begin{tabular}{|c|c|}\hline
  \multicolumn{1}{|c|}{$H=3.2$}   & \multicolumn{1}{|c|}{$H=12.8$} \\ \hline
    \multicolumn{1}{|c|}{Without generalized splines, constant scaling}   & \multicolumn{1}{|c|}{With generalized splines, constant scaling} \\ \hline
     \includegraphics[width=0.45\linewidth]{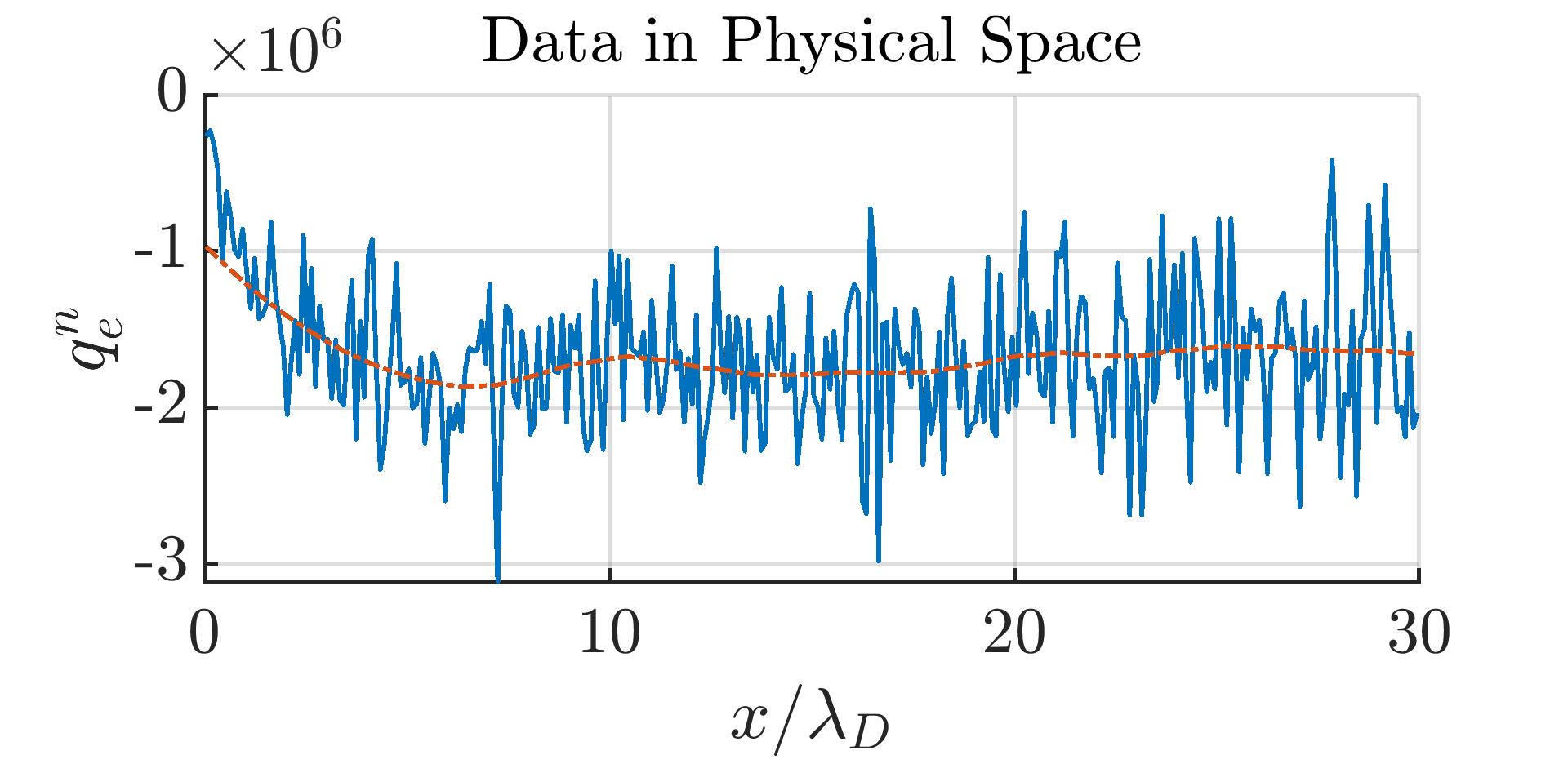} &  \includegraphics[width=0.45\linewidth]{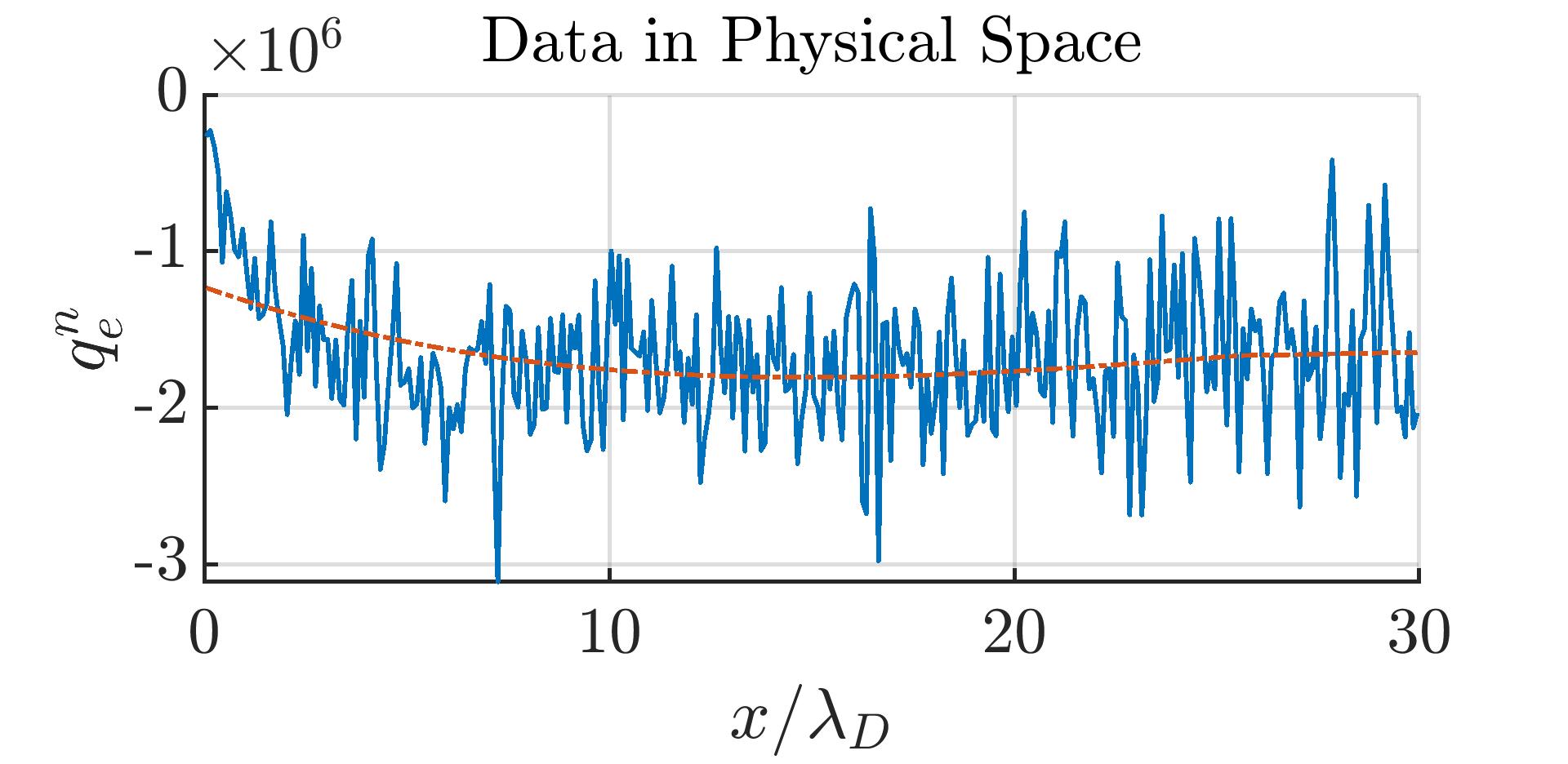}\\ \hline
         \multicolumn{1}{|c|}{With generalized splines, constant scaling}   & \multicolumn{1}{|c|}{With generalized splines and adaptive scaling} \\ \hline
      \includegraphics[width=0.45\linewidth]{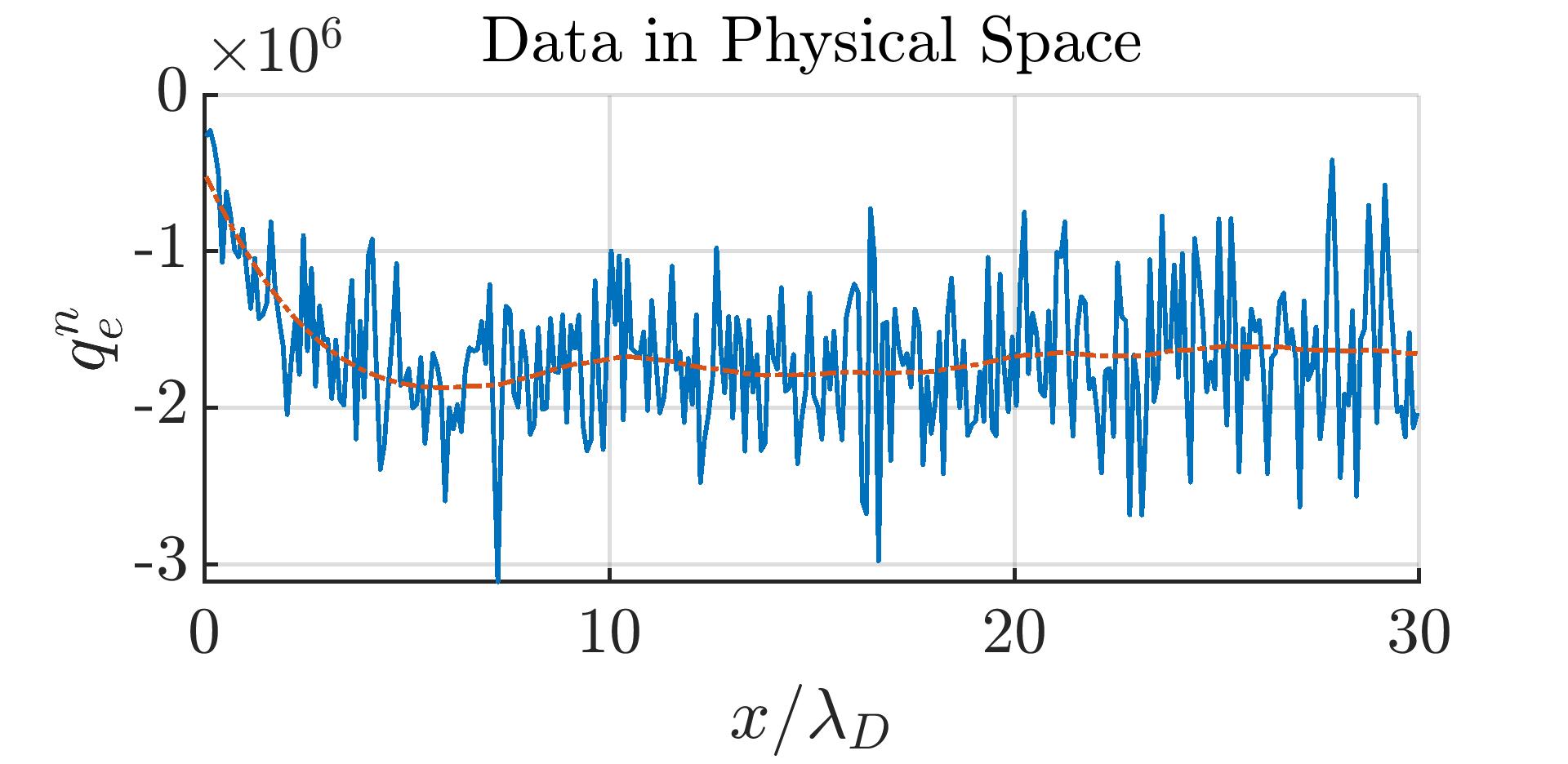} &  \includegraphics[width=0.45\linewidth]{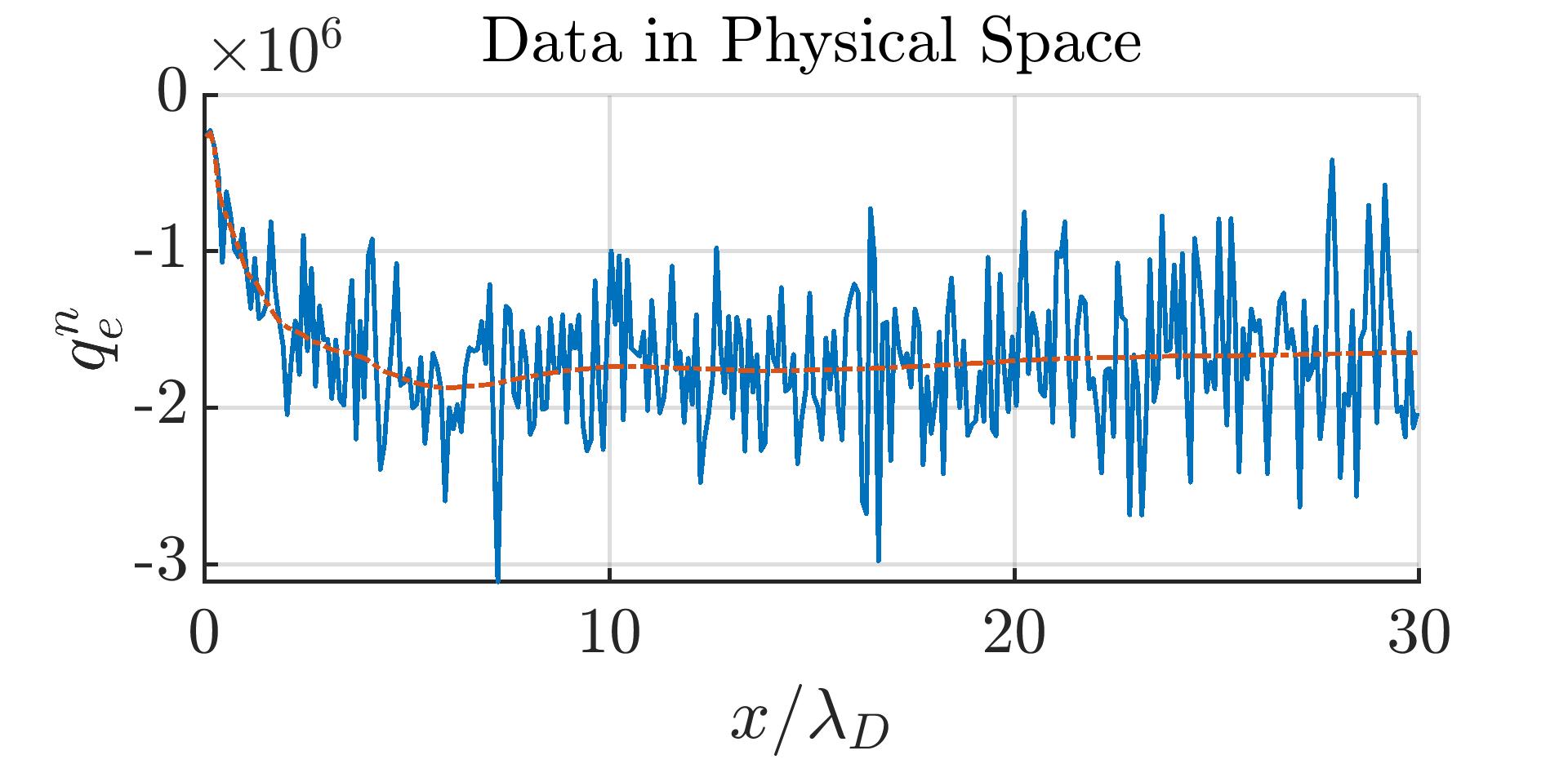}\\ \hline
\end{tabular}
    \centering
    \caption{\label{fig:bnd_behavior_q_n_e}Effect of generalized splines and adaptive kernel scaling on boundary behavior preservation. The underlying data is electron heat flux data $q^n_e$ with $\Omega=[0,800]$ and $h_{grid}=\Delta x=0.1$. Here the solid blue lines denote the unfiltered data and the dashed red lines the filtered data.}
\end{figure}

The second test considers the capability of SIAC filtering to handle large gradients near the boundary. 
This can be a critical characteristic in some cases. Here the considered case is that of the non-neutral sheath (which forms near the region where a plasma intercepts a solid surface) that introduces sharp gradients to the plasma profiles;
for instance see  Fig.~\ref{fig:bnd_behavior_ne} and Fig.~\ref{fig:bnd_behavior_q_n_e}. 
It is worth noting that the electron
conduction heat flux, $q_{n}^e$, in the form of its spatial gradient as shown in Eq.~\eqref{eq-beta}, plays an
important role in setting the Bohm speed~\cite{li-Bohm-PRL,tang2016}, denoising of which is extremely demanding
 since itself suffers higher noise pollution
than the plasma density and flow. 

 Here,  position-dependent SIAC filters that  include generalized splines
and adaptive kernel scalings {are used} to denoise the data so that the Bohm
criterion can be calculated. This avoids a long time average {as} proposed in Ref.~\cite{li-Bohm-PRL},
significantly reducing the computational cost.

The boundary behavior of the filtered data for varying kernel treatments is first examined in the physical space, shown in Fig.~\ref{fig:bnd_behavior_ne}. It can be seen that the typical position-dependent filter composed of 3 B-splines of order 2 with a scaling of $H=32 h_{grid}$ cannot preserve the boundary behavior of electron density. However, if a generalized spline {is added} to the position-dependent kernel to lend more weight to the filtering point in the weighted-average, it is possible to better maintain the boundary behavior. However, if we keep increasing the kernel support to $H=128h_{grid}$ while keeping all the other parameters the same, it can be observed that the generalized spline alone is not enough to maintain {the} boundary behavior. 
This is resolved by using the proposed adaptive kernel scaling given as Eq.~\eqref{eqn:adapt_scale}, as confirmed in Fig.~\ref{fig:bnd_behavior_ne}.

\begin{figure}[tp!]
\begin{tabular}{|c|c|}\hline
  \multicolumn{1}{|c|}{$H=3.2$}   & \multicolumn{1}{|c|}{$H=12.8$} \\ \hline
    \multicolumn{1}{|c|}{Without generalized splines, constant scaling}   & \multicolumn{1}{|c|}{With generalized splines, constant scaling} \\ \hline
     \includegraphics[width=0.45\linewidth]{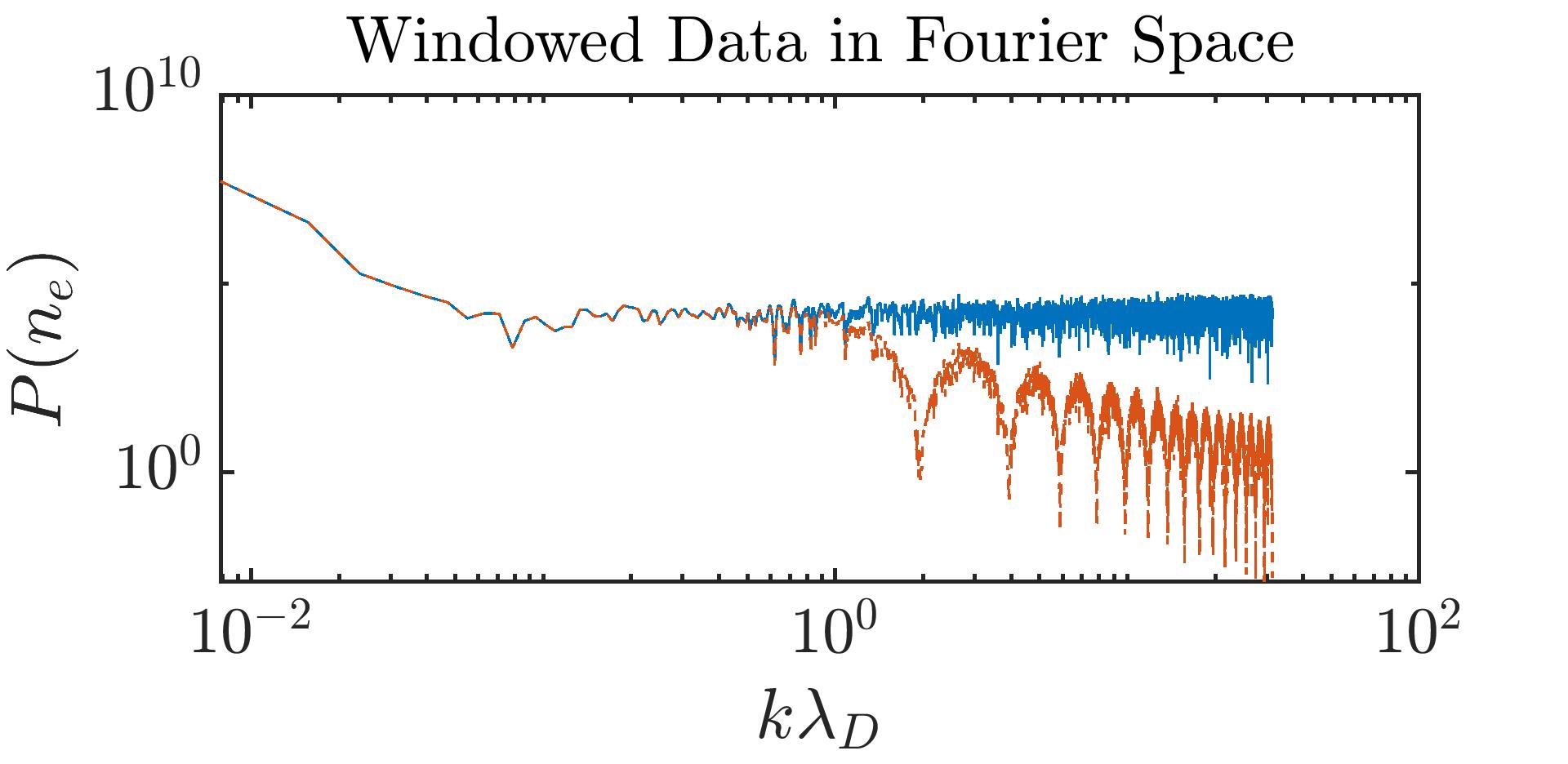} &  \includegraphics[width=0.45\linewidth]{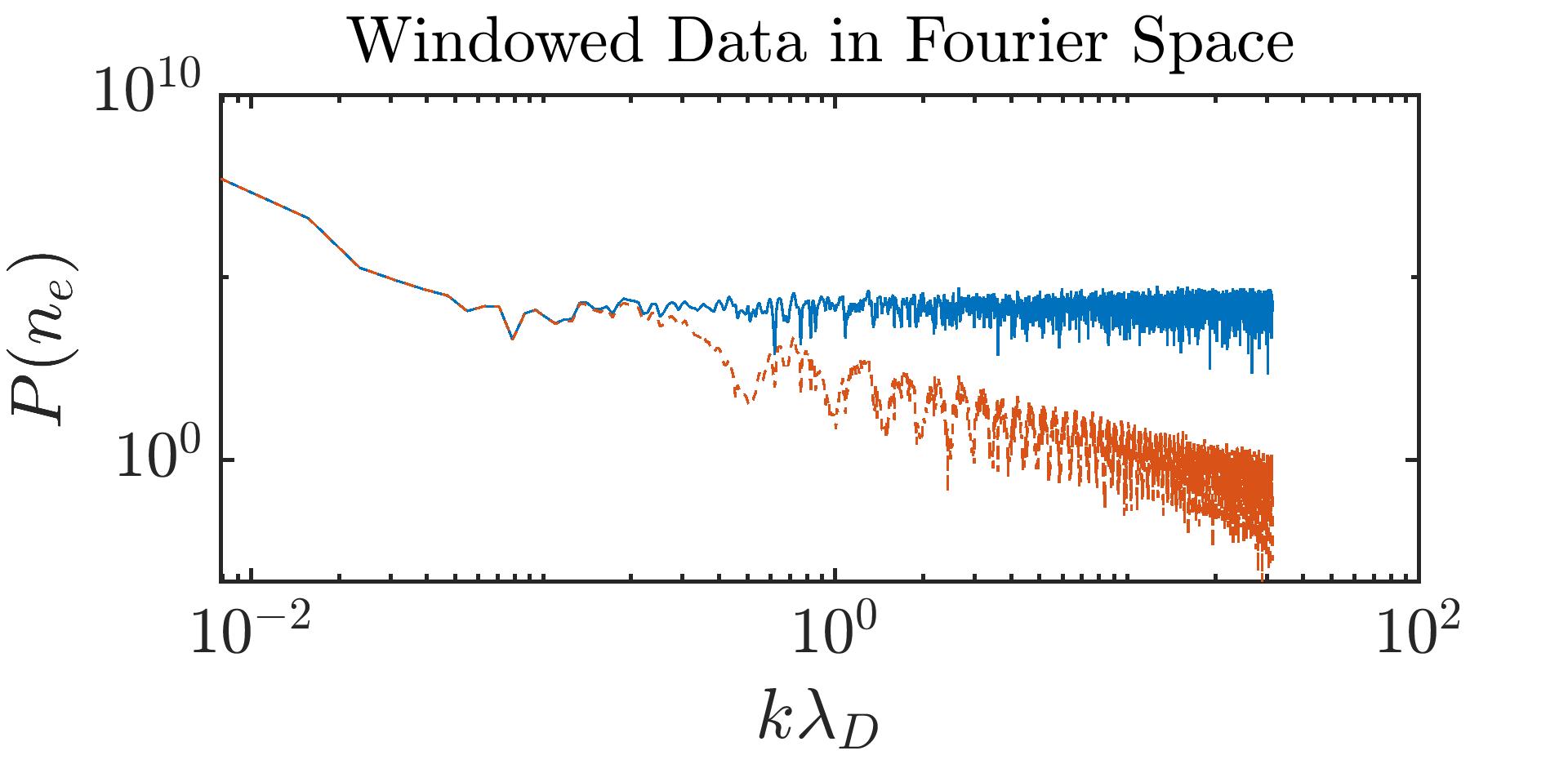}\\ \hline
         \multicolumn{1}{|c|}{With generalized splines, constant scaling}   & \multicolumn{1}{|c|}{With generalized splines and adaptive scaling} \\ \hline
      \includegraphics[width=0.45\linewidth]{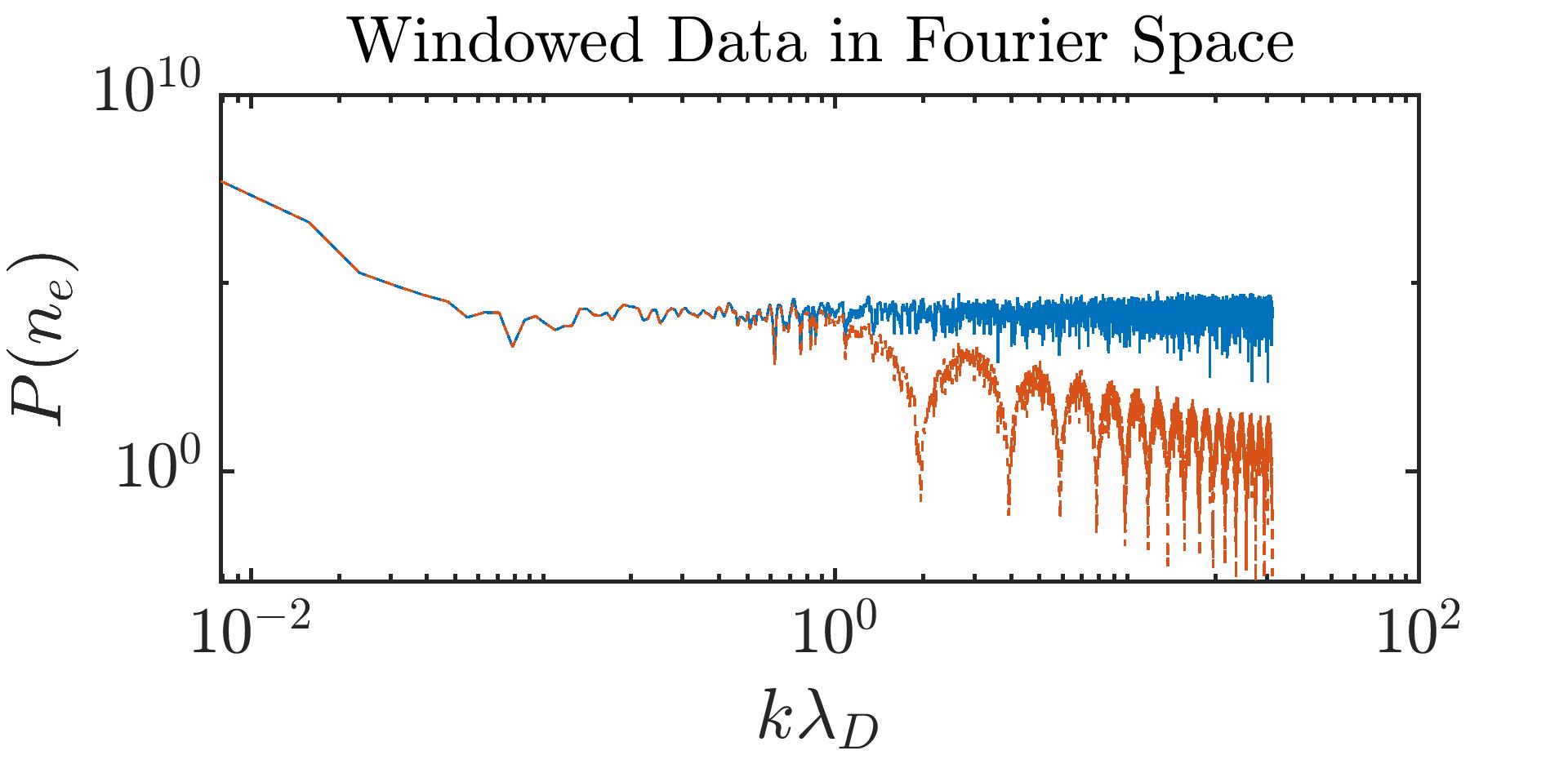} &  \includegraphics[width=0.45\linewidth]{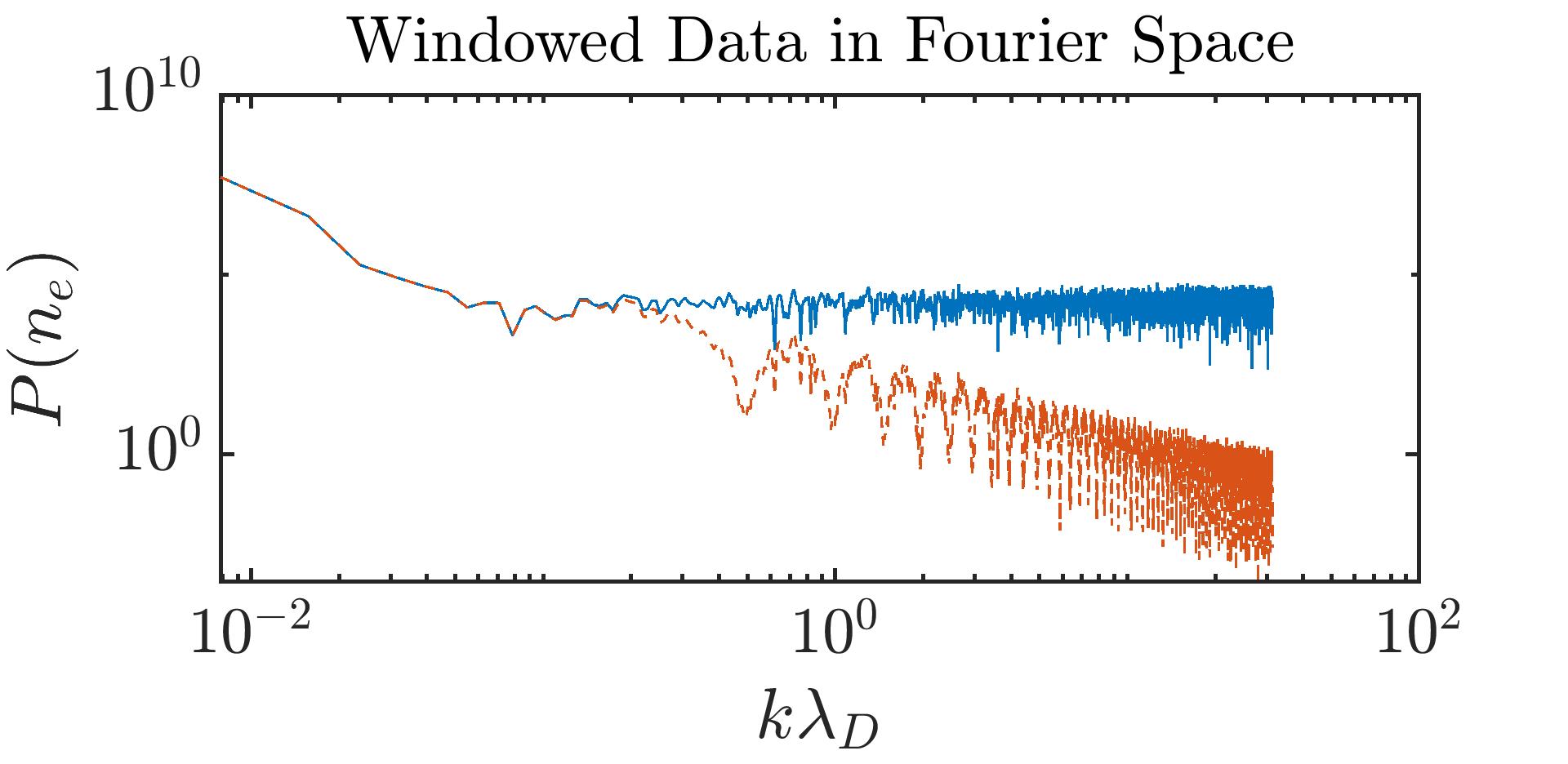}\\ \hline
\end{tabular}
    \centering
    \caption{\label{fig:freq_behavior_ne} Effect of generalized splines and adaptive kernel scaling on the single-sided amplitude spectrum of electron density data $n_e$ with $\Omega=[0,800]$ and $h_{grid}=\Delta x=0.1$. Here the solid blue lines denote the unfiltered data and the dashed red lines the filtered data. Note that the underlying data was non-periodic so a Hanning window function was applied.}
\end{figure}

We then compare the single-sided amplitude spectrums of the filtered and unfiltered electron density data,  results of which are shown in Fig.~\ref{fig:freq_behavior_ne}. 
While larger scalings are shown to have greater damping capabilities,  there is not much of a disparity between the different boundary treatments using the same scaling. This makes sense as, while the differences in tracking the boundary behavior is apparent in physical space, only the data in a small portion of the domain near the boundary varies between these treatments. 
In aggregate, neither adaptive scalings nor generalized splines seem to matter as much as
the choice of interior kernel scaling, when the spectrum of the filtered data is considered. 
Thus, the choice of boundary treatment may not be apparent in spectral analysis.
The results in the Fourier space is largely comparable to the periodic case. 

\begin{figure}[tp!]
\begin{tabular}{|c|c|}\hline
  \multicolumn{1}{|c|}{$H=3.2$}   & \multicolumn{1}{|c|}{$H=12.8$} \\ \hline
    \multicolumn{1}{|c|}{Without generalized splines, constant scaling}   & \multicolumn{1}{|c|}{With generalized splines, constant scaling} \\ \hline
     \includegraphics[width=0.45\linewidth]{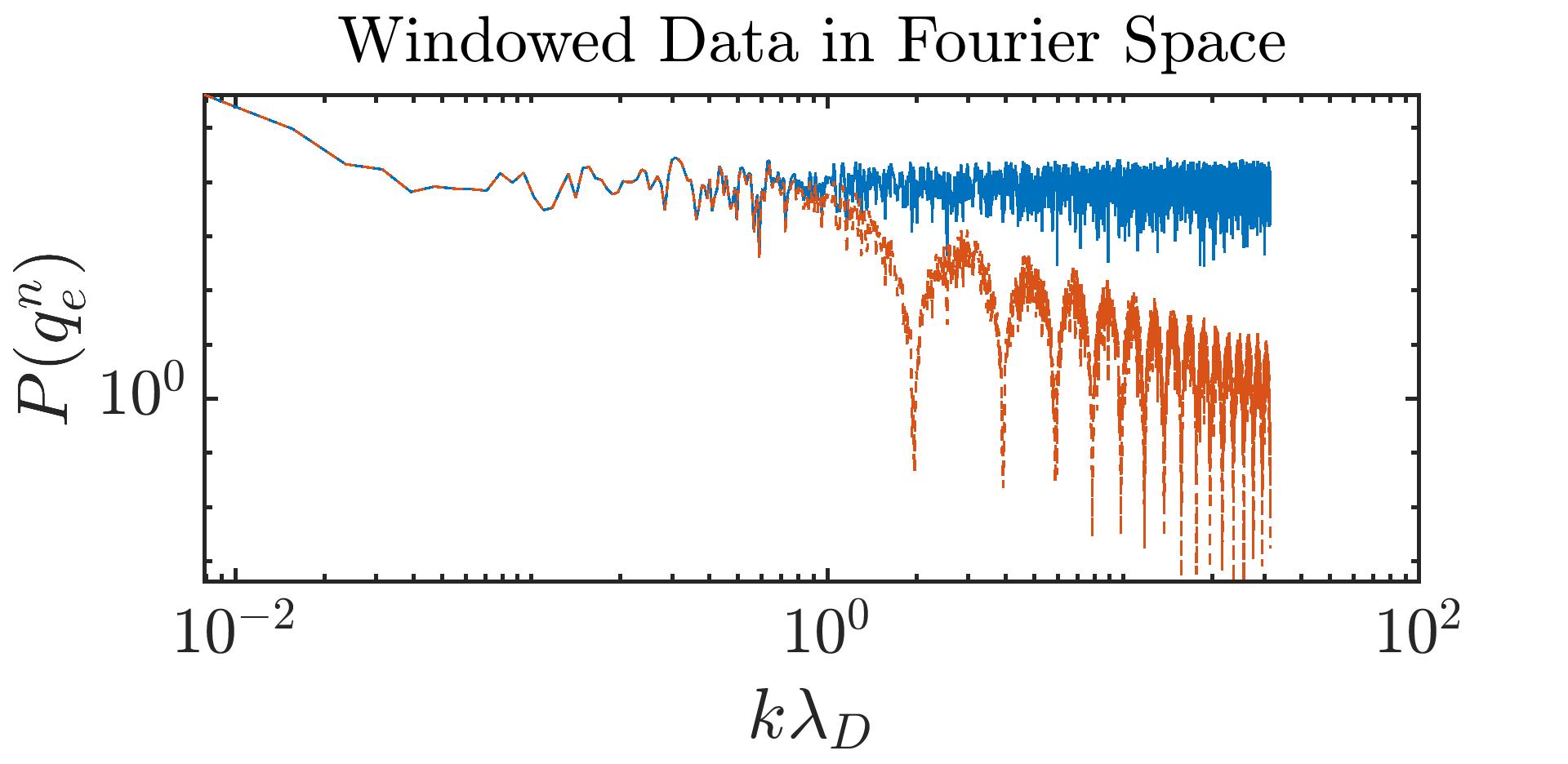} &  \includegraphics[width=0.45\linewidth]{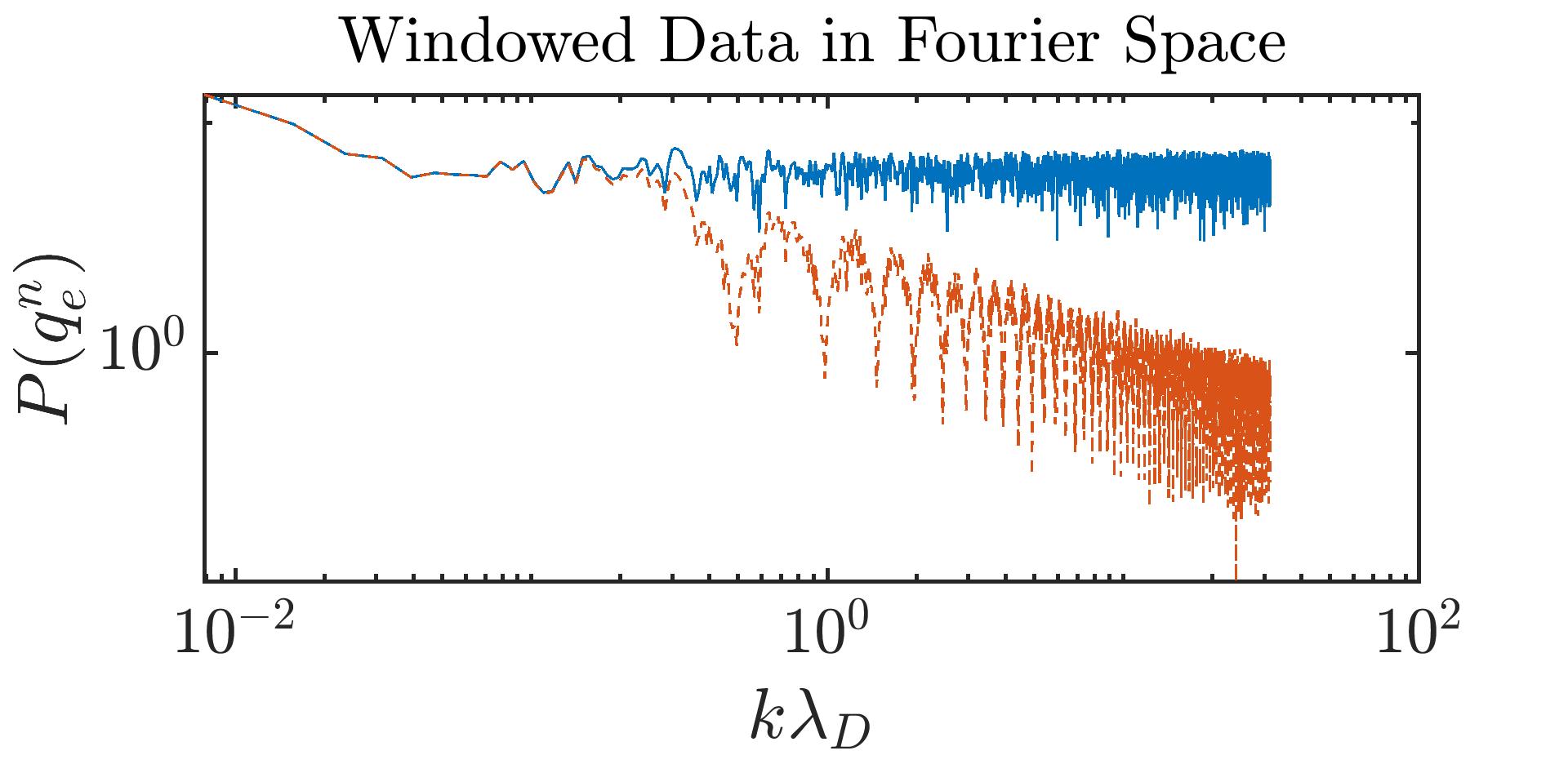}\\ \hline
         \multicolumn{1}{|c|}{With generalized splines, constant scaling}   & \multicolumn{1}{|c|}{With generalized splines and adaptive scaling} \\ \hline
      \includegraphics[width=0.45\linewidth]{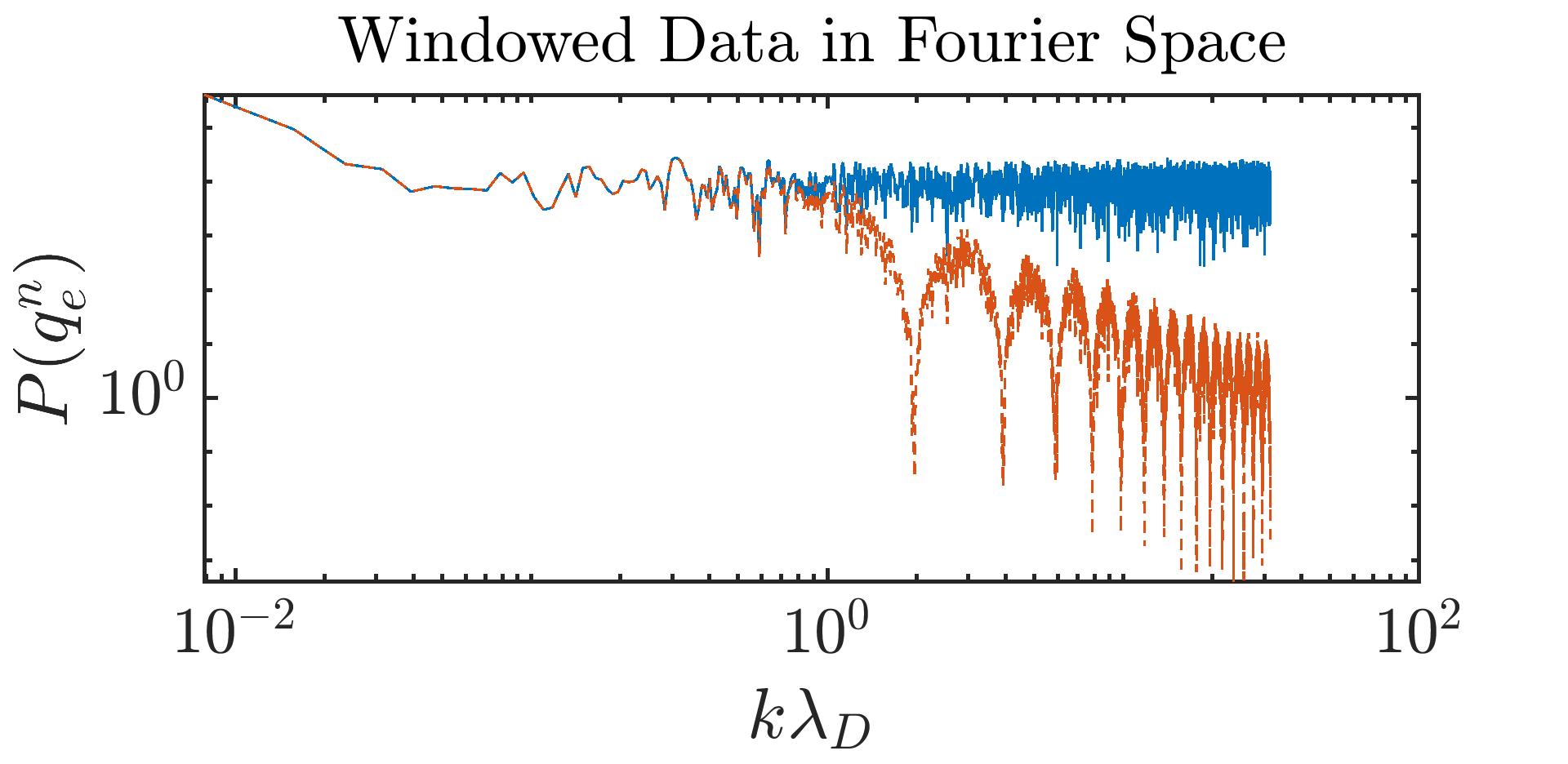} &  \includegraphics[width=0.45\linewidth]{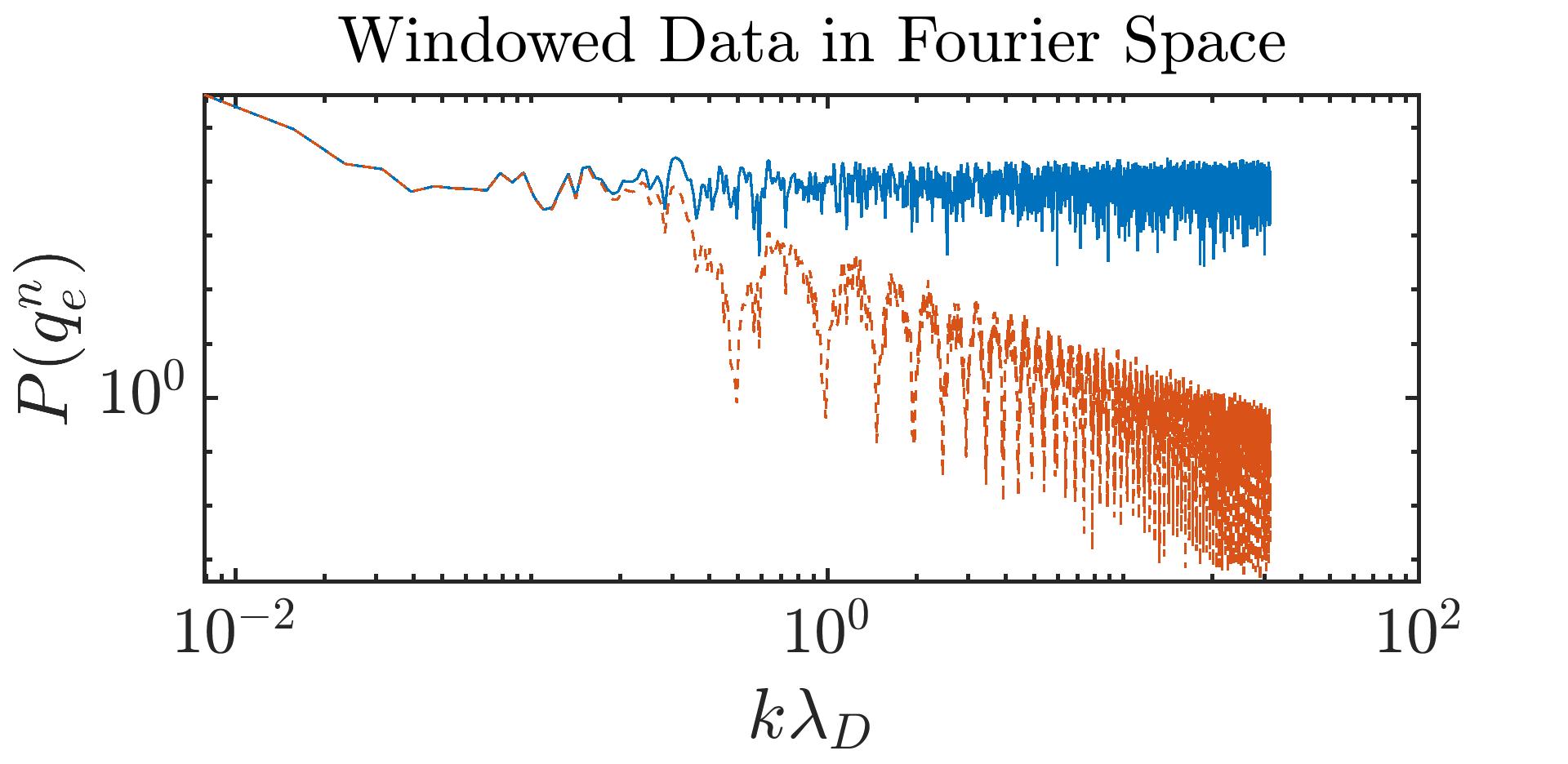}\\ \hline
\end{tabular}
    \centering
    \caption{ \label{fig:freq_behavior_q_n_e}Effect of generalized splines and adaptive kernel scaling on the single-sided amplitude spectrum of electron heat flux data $q^n_e$ with $\Omega=[0,800]$ and $h_{grid}=\Delta x=0.1$. Here the solid blue lines denote the unfiltered data and the dashed red lines the filtered data. Note that the underlying data was non-periodic so a Hanning window function was applied.}
\end{figure}

\begin{figure}[tp!]
\includegraphics[width=0.9\linewidth]{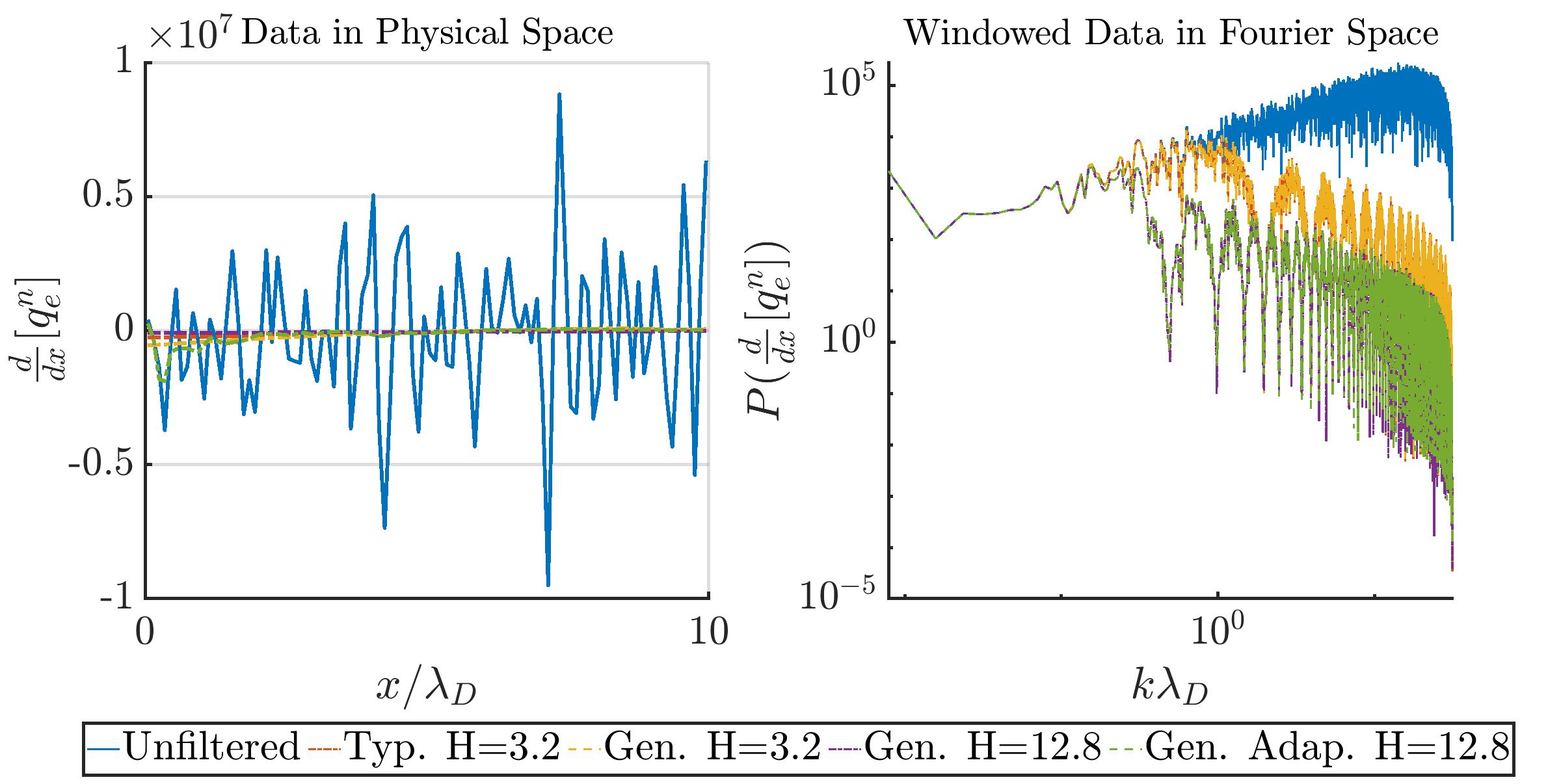}
    \centering
    \caption{\label{fig:ddx_q}Effect of generalized splines and adaptive kernel scaling on electron heat flux gradient data $dq^n_e/dx$ in physical and frequency space. Here $\Omega=[0,800]$, $h_{grid}=\Delta x=0.1$, and the dashed red line (Typ.) denotes a filtering procedure with the typical kernel without generalized splines, while (Gen.) denotes a filtering with a kernel including a generalized spline. Finally, (Adap.) denotes that an adaptive kernel scaling is employed. In physical space, the effects of the differing boundary treatments are most apparent closest to the boundary, so the interior of domain is not shown. Note that the underlying data was non-periodic so a Hanning window function was applied for the spectrum computation.}
\end{figure}

Next  an investigation of the same choice of kernel parameters on the electron heat flux, $q^n_e$, and its spatial gradient, $dq^n_e/dx$, is performed. Note that here in the unfiltered case $dq^n_e/dx$ is computed by a central difference stencil of $q^n_e$ in the interior of the domain, and a 1st order biased stencil at the domain boundaries. The filtered derivative values are obtained by applying the same differencing procedure to the filtered $q^n_e$ values. 
Fig.~\ref{fig:bnd_behavior_q_n_e} and Fig.~\ref{fig:ddx_q} present the effect of the filtering scheme on $q^n_e$ and $dq^n_e/dx$ {is shown}. 
The data is found to be much more oscillatory than the electron density data, though the same behavior {can be seen} for $q^n_e$ under different filter parameters. 
Specifically, adding a generalized spline helps maintain boundary behavior for smaller $H$, but an adaptive scaling becomes neccessary for larger interior scalings. Similar behavior to the $n_e$ case is also observed in the single-sided amplitude spectra depicted in Fig.~\ref{fig:freq_behavior_q_n_e}.  With respect to the gradient of the heat flux, the differing boundary treatments are much less apparent in physical space (see Fig.~\ref{fig:ddx_q}). The presense of noise is particularly apparent in $dq^n_e/dx$ in that it oscillates rapidly about $0$, and much of the energy in the single-sided amplitude spectrum is contained in this noise manifesting in higher frequencies as depicted in Figure \ref{fig:ddx_q}. The application of the filtering procedure serves to damp these higher frequencies. While it may appear that the choice of boundary treatment has little effect on $dq^n_e/dx$, in the next test{, it} will {be seen} that minor perturbations of this variable can have a large impact on computed quantities such as the Bohm speed.

 This result suggests that the SIAC filters can be applied to the non-periodic data, capturing the information well in both the physical space (e.g., boundary layers) and Fourier space. 
 It is particularly encouraging to observe the success in the Fourier space for an extremely noisy case of  $dq^n_e/dx$.


\subsection{Application to Bohm speed calculations}
The numerical verification of the Bohm speed expression of Eq.~\eqref{Bohm_eqn} is a very challenging task. For instance,
the Bohm speed is very sensitive to small perturbations in that if small perturbations like PIC noises in the data cause the denominator of Eq.~\eqref{eq-beta} to flip signs or be near $0$, the computed Bohm speed can be either complex or blow-up.
The contribution to such sensitivity of the Bohm speed to the noises are highly biased in the variables as seen from Eq.~\eqref{eq-beta}, where the electron heat flux gradient, $q_{n}^e$,  contributes the most. On the other hand, $q_{n}^e$ itself  has most of the noise due to it being the high moment calculated from the electron distribution. 
To mitigate such a challenge, Ref.~\cite{li-Bohm-PRL} evolved to a static state solution in an excessively long time, collecting 200,000 snapshots of the static state solution. 
Their average is eventually used to verify the physical law sucessfully. 
Here we aim to show such an expensive brute-force strategy can be avoided by using a small amount of snapshots with the proposed SIAC filters. 
 
 \begin{table}[ht]
    \centering
    \begin{tabular}{|c|c|c|c|c|c|c|c|c|c|c|c|c|c|} \hline
       Variable  &         $q_{n}^e$       &  $q_{n}^i$ &   $Q_{ee}$  &      $Q_{ei}$      &   $Q_{ii}$  &   $E$   &    {$R_T$ }  &     $n_e$ &     $n_i$ &
                $u_{ex}$ &
                $u_{ix}$  &
                 $T_{ex}$ & 
                 $T_{ix}$  \\ \hline
$H$ & 8 &4&16&16&16&2&32&6&6&8&8&16&4\\ \hline
    \end{tabular}
    \caption{\label{tab:scaling_ranges}Variable specific scalings for 2,600 snapshots averaged data set found to produce qualitative agreement with respect to the long-time averaged data.}
\end{table}
 
First, it is necessary to systematically examine the filtering requirements and choose proper scalings for each variable to retain the {majority of the} physics.
Our idea behind these individualized variable scalings is to filter the short-time averaged variables, obtaining filtered data which qualitatively resembles the long-time averaged variables. Procedurally, this takes the form of replacing a single long-time averaged variable in the Bohm speed computation with a noisier short-time averaged counterpart, and then varying the kernel scaling of that short-time averaged variable until good qualitative agreement is obtained with respect to the long-time averaged variable and the long-time averaged Bohm speed. The number of snapshots, 2,600, is the minimum number of time steps for which this can be achieved for all variables. 
The discovered kernel scalings for each component is summarized in Tab.~\ref{tab:scaling_ranges}.
Choosing the scalings in this manner allows for the preservation of the qualitative behavior of the data and real-valued Bohm speed calculations. 

\begin{figure}[ht]
\centering
\begin{tabular}{|c|}\hline
\multicolumn{1}{|c|}{Unfiltered short-time averaged data}\\ \hline
       \includegraphics[width=0.65\linewidth]{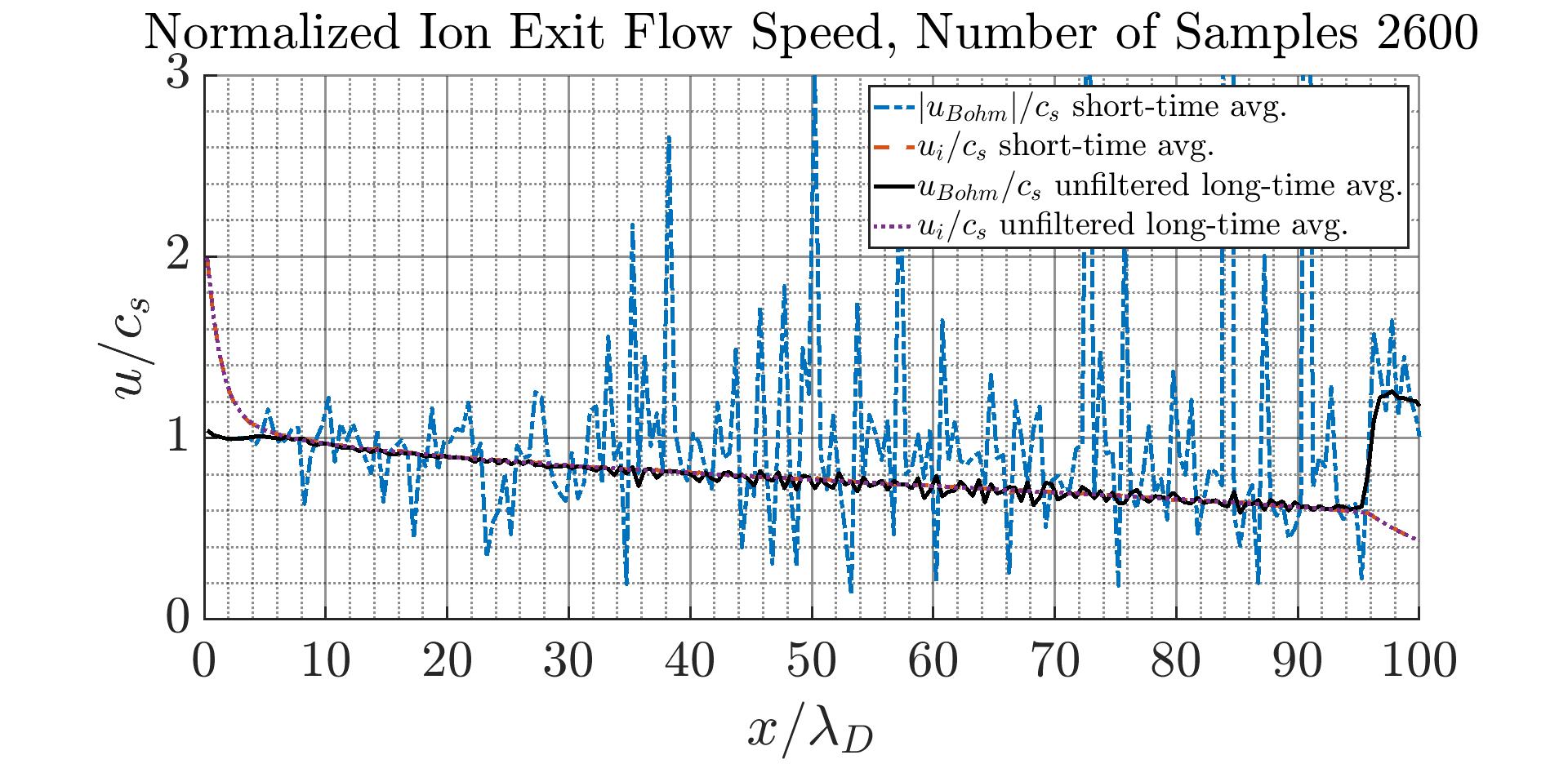}\\ \hline
\multicolumn{1}{|c|}{Filtered short-time averaged data}\\ \hline
         \includegraphics[width=0.65\linewidth]{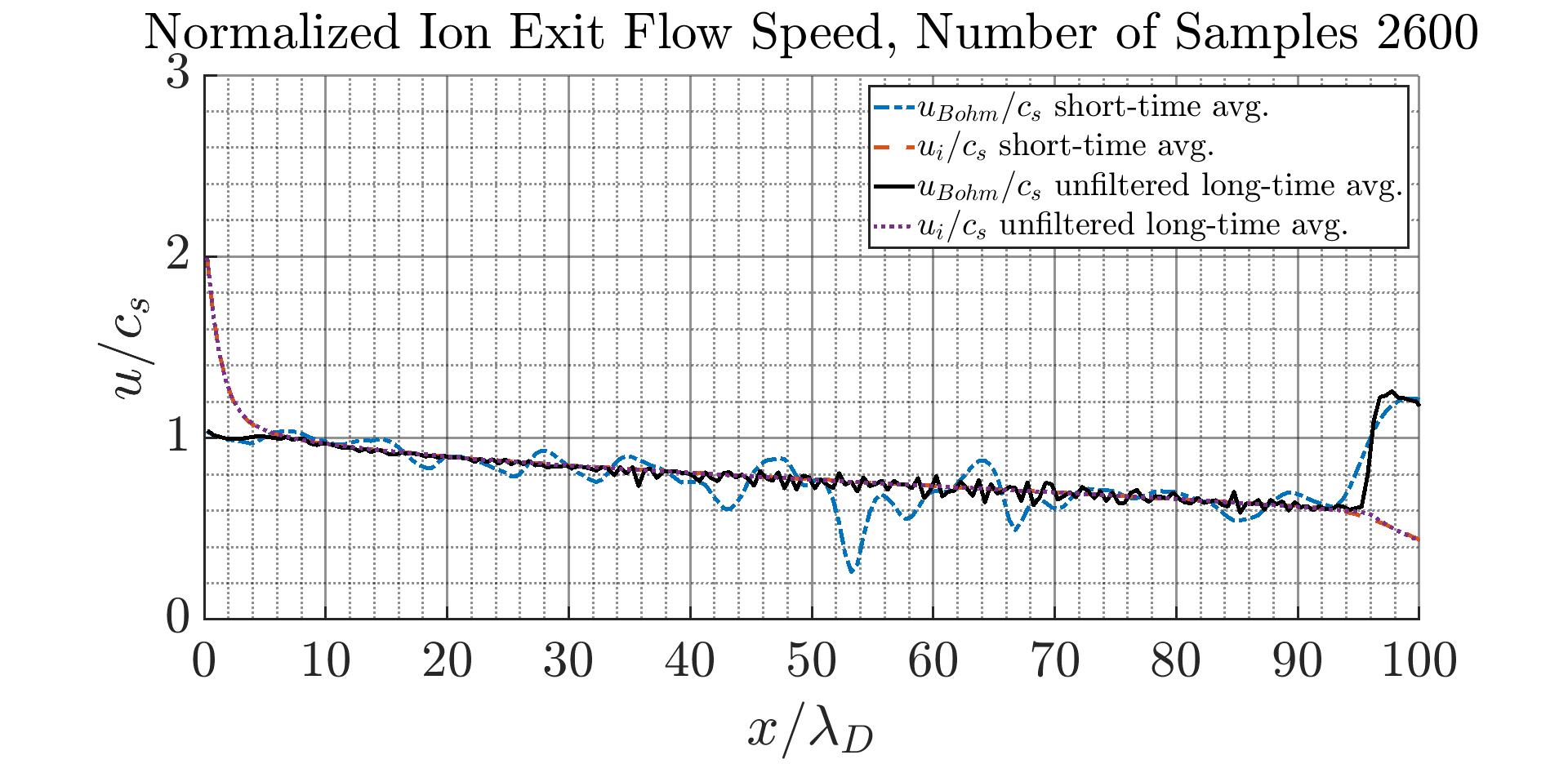}\\ \hline
    \end{tabular}
    \caption{\label{fig:Bohm_comp}Comparison of Bohm speed and ion exit flow using filtered and unfiltered data time-averaged over 2,600 snapshots. The unfiltered long-time averaged data over 
    200,000 snapshots is also plotted. A $K^{(3,2)}_H$ kernel with scalings as indicated in Table \ref{tab:scaling_ranges} was used. Note that the modulus of $u_{Bohm}$ is plotted for the unfiltered short-time average data as the noise was sufficient to induce complex values.}
\end{figure}

Finally, when short-time averaged filtered data alone {is used} with the scalings chosen as indicated in Tab.~\ref{tab:scaling_ranges}, it is capable to obtain a desirable Bohm speed that tracks well with the long-time averaged calculations. This is shown in Fig.~\ref{fig:Bohm_comp}, where  a comparison of the computed Bohm speed from the filtered and unfiltered data with all variables originating from a time average over 2,600 snapshots {is shown}. 
The unfiltered long-time averaged data over 200,000 snapshots is also plotted as a baseline. 
The lower plot of Fig.~\ref{fig:Bohm_comp} shows that the filtered average of normalized Bohm speed (blue dash line) 
matches well with the long-time averaged value (solid black line) from~\cite{li-Bohm-PRL}. 

This result shows that the proposed SIAC filters can be a promising data processing technique in a very practical PIC setting,
which can potentially save significant computational cost. 
For instance, the presented study of Bohm speed suggests that the previous expensive 1D3V PIC simulation in~\cite{li-Bohm-PRL} can be performed in a much shorter time (1.3\%) to discover the same physical law.

  \section{Conclusions and Future Work}

In this paper,  the Smoothness-Increasing Accuracy-Conserving filters have been first-ever extended to effectively denoising  discrete {PIC} data while preserving the underlying physics. 
Several aspects of the filters related to both the periodic boundary condition and  a physical boundary condition that introduces the sharp gradient of plasma profiles have been investigated.
Leveraging knowledge from SIAC applications in a FEM context, the position-dependent filters have been repurposed and  a new adaptive scaling methodology has been developed. 
The proposed filter is assessed using the challenging task of uncovering physical laws, specifically the Bohm speed formulation. 
The outcomes show that the computational cost in PIC can be substantially reduced by reducing the macro-particle numbers or the time steps used for time-averaging.
This underscores the substantial potential of the proposed filters as an effective tool for processing PIC data.

In the future, a more in-depth investigation of the sensitivity of the Bohm speed calculation to noisy variables will be investigated. One more vital capability of the filter that needs to be developed is to self-detect and resolve the interior structures possessing sharp gradients like the shocks, which is of great importance since the time-average approach does not apply to such a dynamical system. Additionally, quantifying the Fourier effects of the adaptively scaled and position-dependent filters will be done.

\section{Acknowledgements}
The first author would like to thank the LANL LDRD ISTI Student fellow program 
and, {along with the fourth author,} the AFOSR under grant number FA9550-20-1-0166 which enabled this work, 
as well as Dr.~Ayaboe Edoh for his insights. 
This work was partially supported by the U.S. Department of Energy SciDAC partnership on Tokamak Disruption Simulation between the
Programs of Fusion Energy Sciences (FES) and Advanced Scientific Computing Research (ASCR). It was also partially supported
by Mathematical Multifaceted Integrated Capability Center (MMICC) of ASCR.
Los Alamos National Laboratory is operated by Triad National Security, LLC, for the National Nuclear Security
Administration of U.S. Department of Energy (Contract No. 89233218CNA000001).
The authors wish to acknowledge and thank Dr.~Yuzhi Li for providing the PIC data on which the Bohm speed calculations were based.

\bibliographystyle{model1_num_names}
\bibliography{main_clean}

\end{document}